\documentclass[12pt,letter]{amsart}
\usepackage[applemac]{inputenc}
\usepackage{amsmath,amsthm,amssymb,latexsym}
\usepackage[T1]{fontenc}
\usepackage{libertine} 
\usepackage{euler} 
\usepackage{eucal} 
\usepackage[pctex32]{color}
\usepackage[pctex32]{graphicx}
\usepackage{pb-diagram}
\usepackage{pst-3d,pst-coil,pst-eps,pst-fill,pst-grad,pst-math,pst-text,pst-tree}

\newcommand\C[1]{\mathcal{#1}}

\newtheorem{proposition}{Proposition}[section]

\newtheorem{assumption}[proposition]{Assumption}
\newtheorem{lemma}[proposition]{Lemma}
\newtheorem{corollary}[proposition]{Corollary}
\newtheorem{fact}[proposition]{Fact}
\newtheorem{claim}[proposition]{Claim}
\newtheorem{question}[proposition]{Question}
\theoremstyle{definition}
\newtheorem{definition}[proposition]{Definition}
\newtheorem{example}[proposition]{Example}
\newtheorem{examples}[proposition]{Examples}
\newtheorem{remark}[proposition]{Remark}
\newtheorem{notation}[proposition]{Notation}

\def\monster{\mathbb{M}}
\def\E{\varepsilon}
\def\dist{\mathbf{d}}
\newcommand{\indep}{\makebox[13pt]{\raisebox{-0.75ex}{$\smile$}\hspace{-3.50mm}
$\mid$\, }}
\newcommand{\eindep}{\makebox[13pt]{\raisebox{-0.75ex}{$\smile$}\hspace{-3.50mm}
$\mid$\, }^{\hspace{-1mm}\varepsilon}}

\def\ordencea{\prec_{\C{K}}}
\def\ordenceac{\succ_{\C{K}}}
\newcommand{\eop}[1]{
\hspace{10mm} \vspace{-6mm}
\begin{flushright}
\qedsymbol$_{\text{#1}}$\\ \ \\
\end{flushright}
}
\newenvironment{prueba}[1][{\it Proof}]{\noindent {\it #1.} }{}
\newcommand{\bdem}[1][Proof]{\begin{prueba}[#1]}
\newcommand{\edem}[1][]{\eop{#1}
\end{prueba}}
\def\bsdem{\begin{prueba}[Reference]}\def\bsindem{\begin{proof}[\ ]}

\def\tuple{\overline}
\def\rest{\upharpoonright}

\def\dominates{\vartriangleright}

\mathchardef\mhyphen="2D
\newcommand{\gatp}{\sf{ga\mhyphen tp}}
\newcommand{\gaS}{\sf{ga\mhyphen S}}

\def\dominates{\vartriangleright}
\def\isdominated{\vartriangleleft}
\def\equidominated{\dominates\hspace{-1mm}\isdominated}
\newcommand{\stype}{\mathfrak{S}t}

\begin{document}
\title[Around Independence and Domination in MAEC]{Around Independence
  and Domination in Metric Abstract Elementary Classes,
under Uniqueness of Limit Models.}
\author[A. Villaveces and P. Zambrano]{Andr\'es Villaveces\\ Pedro Zambrano}
\address{{\rm E-mail - A. Villaveces:} {\it avillavecesn@unal.edu.co}\\
Departamento de Matem\'aticas Universidad Nacional de Colombia, AK
30 $\#$ 45-03, Bogot\'a - Colombia}

\address{{\rm E-mail - P. Zambrano:} {\it phzambranor@unal.edu.co}\\
Departamento de Matem\'aticas Universidad Nacional de Colombia, AK
30 $\#$ 45-03, Bogot\'a - Colombia}
\thanks{The second author wants
  to thank the first author for the time devoted to advising him for
  his Ph.D.'s thesis, one of whose fruits is this paper. He also
  thanks Tapani Hyttinen for helpful discussions in Helsinki
  during November 2008 and August-October 2009 and thanks John Baldwin
  for the nice discussions in Bogotá on domination in (discrete) AECs
  in November 2009. Both authors were partially supported by
  Colciencias.}
\date{\today}
\begin{abstract}
We study notions of independence appropriate for a stability theory of
metric abstract elementary classes (for short, MAECs). We build on
previous notions used
in the discrete case, and adapt definitions to the metric case. In
particular, we study notions that behave well under
superstability-like assumptions. Also, under uniqueness of limit
models, we study domination, orthogonality and parallelism of Galois
types in MAECs.
\end{abstract}
\maketitle

\section*{Introduction}
\label{sec:intro}

In the study of the Stability Theory of Abstract Elementary Classes
(for short, AEC, in this paper),
various versions of independence linked to
\emph{splitting} (introduced originally by Shelah in the discrete AEC
case~\cite{Sh 394}) have played an important r\^{o}le. Various
categoricity transfer results, as well as the development of stability
theory in AEC have so far used non-splitting as one of the main independence
notions.

In the metric continuous case (a generalization of both
usual, or ``discrete'' AEC \emph{and} ``First Order''
Continuous Model Theory), notions of independence have been used with
some success in a strongly homogeneous $\omega$-stable,
(L\"{o}wenheim-Skolem number $\aleph_0$) case by \r{A}sa
Hirvonen~\cite{Hi}.

We focus here in a notion of independence, called
\emph{smooth independence} (see~\ref{r_independence}), that
generalizes non-splitting
to the metric context,
 and
works well under the existence of various sorts of limit models. We
study conditions under which smooth independence satisfies appropriate
variants of transitivity,
stationarity, extension, existence (Section~\ref{sec:2}). We also
study the continuity of
this independence notion (see~\ref{cont_indep}).

Applications of these techniques include the study of ``superstable''
metric AEC (limit models and smooth towers~\cite{ViZa}), and steps towards
a generalization of the main theorem in~\cite{GrVaVi}
(Uniqueness of Limit Models in AEC under superstability-like
assumptions under categoricity). They also include notions of domination appropriate for
both discrete and metric continuous (superstable) AEC.
\\ \\
In \cite{Ba}, J. Baldwin does a study of a weak notion of domination which is based on a rough notion of independence in terms of on intersections of models, although he assumes uniqueness of limit models as a superstability-like assumption.
\\ \\
In superstable first order theories, getting a decomposition (up to equi-domination) of stationary types as a finite product of regular types provides us a proof of the following fact due to Lachlan: A countable superstable theory has $1$ or infinitely many countable models. Also, there are versions of this decomposition theorem in not necessarily stable theories (see \cite{OnUs2}) such as rosy and dependent theories (see \cite{OnUs1}), settings where there is a very well-behaved independence notion.
\\ \\
In section 3, we introduce a notion of dominance in the setting of superstable MAEC. We base our work on \cite{Ba}, but we define our notion of dominance using smooth independence and not just using intersections as J. Baldwin does in his paper. We prove that under suitable assumptions, given a tuple $(M,\C{M},\tuple{a},N)$ (where $\C{M}$ is a resolution of $M$ which witnesses that $M$ is a limit model over some model $M_0$) such that $\tuple{a}\indep_{M_0}^{\C{M}_0} M$ (and therefore $\gatp(\tuple{a}/M)$ is a stationary type because $M$ is an universal model over $M_0$), there exist $N^*$ and a resolution $\C{M}^*$ which witnesses that $M$ is a limit model over $M_0$ such that $\tuple{a}\equidominated_{M}^{\C{M}^*}N^*$. Also, in this chapter we study notions of orthogonality and parallelism in superstable MAECs, inspired in \cite{Sh705}. In this study, we drop some strong conditions given in \cite{Sh705} and simplify some of the proofs given there. Also, we prove some properties which were not studied in \cite{Sh705}.


\section{Some basic definitions and results}
\begin{definition}\label{DensityCharacter}
The \emph{density character} of a topological space is the smallest
cardinality of a dense subset of the space. If $X$ is a topological
space, we denote its density character by $dc(X)$. If $A$ is a
subset of a topological space $X$, we define
$dc(A):=dc(\overline{A})$.
\end{definition}

We consider a natural adaptation of the notion of {\it Abstract
Elementary Class} (see \cite{Gr} and \cite{BaMon}), but work in
a context of Continuous Logic that generalizes the ``First Order
Continuous'' setting of~\cite{CoMon} by removing the assumption of
uniform continuity\footnote{Uniform continuity guarantees logical
  compactness in their formalization, but we drop compactness in
  AEC-like settings.}. We follow the definitions given by \r{A}sa Hirvonen
and Tapani Hyttinen (see \cite{Hi}).

\begin{definition}\label{MAEC}
Let $\mathcal{K}$ be a class of $L$-structures (in the context of
Continuous Logic) and $\ordencea$ be a binary relation defined in
$\C{K}$. We say that $(\C{K},\ordencea)$ is a {\it Metric Abstract
Elementary Class} (shortly {\it MAEC}) if:

\begin{enumerate}
\item $\C{K}$ and $\ordencea$ are closed under isomorphism.
\item $\ordencea$ is a partial order in $\C{K}$.
\item If $M\ordencea N$ then $M\subseteq N$.
\item (\emph{Completion of Union of Chains}) If $(M_i:i<\lambda)$ is a $\ordencea$-increasing chain then
    \begin{enumerate}
    \item the function symbols in $L$ can be uniquely interpreted on the completion of
        $\bigcup_{i<\lambda} M_i$ in such a way that
        $\overline{\bigcup_{i<\lambda}M_i} \in \C{K}$
    \item for each $j < \lambda$ , $M_j \ordencea
        \overline{\bigcup_{i<\lambda} M_i}$
    \item if each $M_i\in \C{K}\in N$, then
        $\overline{\bigcup_{i<\lambda} M_i} \ordencea N$.
    \end{enumerate}
    \item (\emph{Coherence}) if $M_1\subseteq M_2\ordencea M_3$ and
        $M_1\ordencea M_3$, then $M_1\ordencea M_2$.
    \item (DLS) There exists a cardinality $LS^d(K)$ (which is called the
      {\it metric L\"{o}\-wen\-heim-Skolem number}) such that if $M
      \in \C{K}$
        and $A \subseteq M$, then there exists $N\in \C{K}$ such that $dc(N) \le dc(A) + LS^d(K)$
        and $A\subseteq N \ordencea M$.
\end{enumerate}
\end{definition}

\begin{examples}
\begin{enumerate}
\item Any continuous elementary class (see~\cite{CoMon}) with the
  usual elementary substructure relation is an MAEC. Important cases
  include
    \begin{enumerate}
        \item Hilbert spaces with a unitary operator (Argoty and
          Berenstein, see~\cite{ArBe}).
        \item Nakano spaces with compact essential rank (Poitevin,
          see~\cite{Po}.
        \item Probability Spaces.
        \item Compact Abstract Theories, see \cite{Be1,Be2}
    \end{enumerate}
\item Gelfand triplets (see \cite{Za} and forthcoming \cite{GaZa}).
\item Hilbert Spaces, with various classes of unbounded operators.
\item A subclass of completions of metric spaces which satisfy
  approximately a positive bounded theory, where $\ordencea$ is
  interpreted by the approximate elementary submodel relation (see
  \cite{HeIo}).
\item Various classes of Banach spaces, where $\ordencea$ is interpreted
by the closed subspace relation\footnote{Notice that
  several of these classes fall under case (1) --- however, in
  general, natural classes of Banach spaces are not axiomatizable in
  the context of~\cite{CoMon}.} (see \cite{Hi}).
\end{enumerate}
\end{examples}

\begin{definition}\label{k_embed}
We call a function $f:M\to N$  a $\mathcal{K}$-\emph{embedding} if
\begin{enumerate}
\item For every $k$-ary function symbol $F$ of $L$, we have $f(F^M(a_1\cdots
  a_k))=F^N(f(a_1)\cdots f(a_k))$.
\item For every constant symbol $c$ of $L$, $f(c^M)=c^N$.
\item For every $m$-ary relation symbol $R$ of $L$, for every
  $\bar{a}\in M^{m}$, $d(\bar{a},R^M)=d(f(\bar{a}),R^N)$.
\item $f[M]\ordencea N$.
\end{enumerate}
\end{definition}

\begin{definition}[Amalgamation Property, AP]\label{AP}
Let $\C{K}$ be an MAEC. We say that $\C{K}$ satisfies {\it Amalgamation Property} (for short {\it AP}) if and only if  for every $M,M_1,M_2\in \C{K}$, if $g_i:M\to M_i$ is a $\C{K}$-embedding $(i\in \{1,2\})$ then there exist $N\in \C{N}$ and $\C{K}$-embeddings $f_i:M_i\to N$ ($i\in \{1,2\}$) such that $f_1\circ g_1=f_2\circ g_2$.
\end{definition}

\[
\begin{diagram}
\node{M_1} \arrow{e,t,..}{f_1} \node{N}\\
\node{M} \arrow{n,l}{g_1} \arrow{e,b}{g_2} \node{M_2} \arrow{n,r,..}{f_2}
\end{diagram}
\]

\begin{definition}[Joint Embedding Property, JEP]
Let $\C{K}$ be an MAEC. We say that $\C{K}$ satisfies {\it Joint Embedding Property} (for short {\it JEP}) if and only if  for every $M_1,M_2\in \C{K}$ there exist $N\in \C{N}$ and $\C{K}$-embeddings $f_i:M_i\to N$ ($i\in \{1,2\}$).
\end{definition}


\begin{remark}
Notice that if $\C{K}$ has a prime model, then AP implies JEP.
\end{remark}

\begin{remark}[Monster Model]\label{Monster_Model}
If $\C{K}$ is an MAEC which satisfies AP and JEP and has large enough models, then we can
construct a large enough model $\monster$ (which we call a {\it
  Monster Model}) which is homogeneous --i.e., every isomorphism
between two $\C{K}$-substructures of $\monster$ can be extended to an
automorphism of $\monster$-- and also universal --i.e., every model
with density character $<dc(\monster)$ can be $\C{K}$-embedded into
$\monster$.
\end{remark}

\begin{definition}[Galois type]\label{Galois_Type}
Under the existence of a monster model $\monster$ as in
remark~\ref{Monster_Model}, for all $\tuple{a}\in \monster$ and
$N\ordencea \monster$, we define $\gatp(\tuple{a}/N)$ (the
\emph{Galois type of $\tuple{a}$ over $N$}) as the orbit of
$\tuple{a}$ under $Aut(\monster/N):=\{f\in Aut(\monster): f\rest N =
id_{N}\}$. We denote the space of Galois types over a model $M\in
\C{K}$ by $\gaS(M)$.
\end{definition}

Throughout this paper, we assume the existence of a homogenous and universal monster model as in remark~\ref{Monster_Model}.

\begin{definition}[Distance between types]\label{Distance_Types}
Let $p,q\in \gaS(M)$. We define $d(p,q):=\inf\{d(\tuple{a},\tuple{b}):
\tuple{a},\tuple{b}\in \monster, \tuple{a}\models p , \tuple{b}\models
q\}$, where ${\rm lg}(\tuple{a})={\rm lg}(\tuple{b})=:n$ and
$d(\tuple{a},\tuple{b}):=\max\{d(a_i,b_i): 1\le i\le n\}$.
\end{definition}

\begin{definition}[Continuity of Types]\label{Continuity_Types}
Let $\C{K}$ be an MAEC and consider $(a_n)\to a$ in $\monster$. We say that $\C{K}$ satisfies {\it Continuity of Types Property}\footnote{This property is also called {\it Perturbation Property in \cite{Hi}}} (for short, {\it CTP}), if and only if, if $\gatp(a_n/M)=\gatp(a_0/M)$ for all $n<\omega$ then $\gatp(a/M)=\gatp(a_0/M)$.
\end{definition}

\begin{remark}
In general, distance between types $d$ (see Definition~\ref{Distance_Types}) is just a pseudo-metric. But it is straightforward to see that the fact that $d$ is a metric is equivalent to CTP.
\end{remark}

Throughout this paper, we also assume CTP (so, distance between types is in fact a metric).

\begin{definition}[Universality]\label{Universality}
Let $\C{K}$ be an MAEC and $N\ordencea M$. We say that $M$ is $\lambda$-{\it universal} over $N$ iff
for every $N'\ordenceac N$ with density character $\lambda$ there exists a $\C{K}$-embedding $f:N'\to M$ such that $f\rest N=id_N$. We say that $M$ is {\it universal} over $N$ if $M$ is $dc(M)$-universal over $N$.
\end{definition}



\begin{lemma}\label{unions}
Let $\C{K}$ be a MAEC. 
If                       $f_i: M_i\to M$ ($i<\mu$) is a
 $\subseteq$-increasing
and continuous (in the metric sense) chain of $\C{K}$-embeddings,
then there exists a $\C{K}$-embedding $f:\overline{\bigcup_{i<\mu}
M_i}\to M$ which extends $g:=\bigcup_{i<\mu}
f_i:\bigcup_{i<\mu}M_i\to M$.
\end{lemma}
\bdem Let $a\in \overline{\bigcup_{i<\mu}M_i}$, so there
exist elements $a_n\in \bigcup_{i<\mu}M_i$ for $n<\omega$, such that
$(a_n)_{n<\omega}\to a$. 
As $(a_n)_{n<\omega}$ is a Cauchy sequence,
$(g(a_n))_{n<\omega}$ is also a Cauchy sequence (since $g$ is an
isometry). So, there exists $b\in M$  such that
$(g(a_n))_{n<\omega}\to b$. Define $f(a):=b$. Proceed in a similar
way for every $a\in \overline{\bigcup_{i<\mu}M_i}$. The function $f$ is
well-defined: if we take $(a'_n)_{n<\omega}$ a sequence in
$\bigcup_{i<\mu}M_i$ such that $(a'_n)_{n<\omega}\to a$, let
$b'\in M$ be such that $(g(a'_n))_{n<\omega}\to b'$. We will
prove that $b=b'$.
Otherwise, let $\varepsilon:=d(b,b')>0$.\\

\begin{claim}\label{an}
Given $\varepsilon'>0$, there exists $N<\omega$ such that for all
$n\ge N$ $d(g(a_n),g(a'_n))<\varepsilon'$.
\end{claim}
\bdem As $(a_n)_{n<\omega}\to a$ and $(a'_n)_{n<\omega}\to a$,
there exists $N<\omega$ such that for all $n\ge N$ we have that
$d(a_n,a)< \varepsilon'/2$ and $d(a'_n,a)<\varepsilon'/2$, so for
all $n\ge N$ we have that $d(a_n,a'_n)\le
d(a_n,a)+d(a,a'_n)<\varepsilon'$. As $g$ is an isometry, for all
$n\ge N$ we have that $d(g(a_n),g(a'_n))<\varepsilon'$.
\edem[Claim~\ref{an}]

As $(g(a_n))_{n<\omega}\to b$, $(g(a'_n))_{n<\omega}\to b'$ and by
claim \ref{an}, there exists $M<\omega$ such that for all $n\ge M$
we have that $d(g(a_n),b)<\varepsilon/3$,
$d(g(a'_n),b')<\varepsilon/3$ and
$d(g(a_n),g(a'_n))<\varepsilon/3$. So, for all $n\ge M$ we have
that $d(b,b')\le
d(b,g(a_n))+d(g(a_n),g(a'_n))+d(g(a'_n),b')<\varepsilon=d(b,b')$
(contradiction).\\ \\
\indent Therefore $b=b'$ and so $f$ is well-defined.
\\ \\
\indent We have that $f$ extends $g$: let $a\in \bigcup_{i<\omega}
M_i $, so taking $a_n:=a$ ($n<\omega$) we have that
$(a_n)_{n<\omega}\to a $ and $(g(a_n))_{n<\omega}$ is also a
constant sequence. So, $f(a):=\lim_{n<\omega} g(a_n)=g(a)$.
\\ \\
 \indent Let $c\in f[\overline{\bigcup_{i<\mu} M_i}]$,
so there exists $a\in \overline{\bigcup_{i<\mu}M_i}$ such that
$f(a)=c$, so there exists $(a_n)_{n<\omega}$ a sequence in
$\bigcup_{i<\mu}M_i$ such that $(a_n)_{n<\omega}\to a$ and
$c:=\lim_{n<\omega} g(a_n)$. Therefore $c\in
\overline{g[\bigcup_{i<\mu}M_i]
}=\overline{\bigcup_{i<\mu}f_i[M_i] }$, so
$f[\overline{\bigcup_{i<\mu}M_i }]\subseteq
\overline{\bigcup_{i<\mu}f_i[M_i] }$. Take $c\in
\overline{\bigcup_{i<\mu}f_i[M_i]}=\overline{g[\bigcup_{i<\mu}M_i]}$,
so there exists a sequence $(b_n)_{n<\omega}$ in
$\bigcup_{i<\mu}M_i$ such that $(g(b_n))_{n<\omega}\to c$. As
$(g(b_n))_{n<\omega}$ is a Cauchy sequence and $g$ is an isometry,
we have that $(b_n)_{n<\omega}$ is also a Cauchy sequence. So,
there exists $a\in \overline{\bigcup_{i<\mu}M_i}$ such that
$(b_n)_{n<\omega}\to a$, and therefore $f(a):=\lim_{n<\omega}
g(b_n)=c $, hence $c\in f[\overline{\bigcup_{i<\mu}M_i} ]$.
So, $f[\overline{\bigcup_{i<\mu}M_i}
]=\overline{\bigcup_{i<\mu}f_i[M_i]}$. As $(f_i:i<\mu)$ is a
$\subseteq$-increasing and continuous chain of $\C{K}$-embeddings,
$f_i[M_i]\ordencea M$, so by completion of union of chains (definition~\ref{MAEC} (4)) and coherence (definition~\ref{MAEC} (5)) MAEC axioms and we have that $f[\overline{\bigcup_{i<\mu}M_i}
]=\overline{\bigcup_{i<\mu}f_i[M_i]}\ordencea M$. Furthermore, for
every symbol $\sigma$ of $L(\mathcal{K})$, $f$ is compatible with the
interpretation of $\sigma$ in $\overline{\bigcup_{i<\mu} M_i}$: $f$ is
a limit of $\mathcal{K}$-embeddings -- function symbols on these
limits are uniquely interpreted by Axiom 4(a), and $f$ being a limit
of $\mathcal{K}$-embeddings, distances to interpretations of
predicates are preserved.
Therefore $f$ is a $\C{K}$-embedding which extends $g$.
\edem[Lemma~\ref{unions}]

\begin{fact}[Hyttinen-Hirvonen]\label{inf}
Given $\varepsilon>0$ and $a\models p$, there exists $b\models q$
such that $d(a,b)\le \dist(p,q)+\varepsilon$
\end{fact}
\bdem Fix $\varepsilon > 0$. By the definition of $d$, there exist
realizations $c\models p$ and $c'\models q$ such that $d(c, c')\le
\dist(p, q) + \varepsilon$. As $a,c\models p $ then there exists
$f\in Aut(\mathbb{M}/A)$ such that $f(c) = a$. Note that
$d(a,f(c'))=d(f(c),f(c'))=d(c,c')\le \dist(p,q)+\varepsilon$,
where $f(c')\models q$, so $f(c')$ is the required $b$.
\edem[Fact~\ref{inf}]

\begin{corollary}\label{inf_consequence}
Given $\E>0$  and $p,q\in \gaS(M)$ such that
$\dist(p,q)<\varepsilon$ and $b\models q$, then there exists
$a_{\varepsilon}\models p$ such that
$d(a_{\varepsilon},b)<2\varepsilon$.
\end{corollary}

\bdem By fact \ref{inf}, there exists $a_{\varepsilon}\models p$
such that $d(a_{\varepsilon},b)\le \dist(p,q)+\varepsilon$,
therefore $d(a_{\varepsilon},b)\le \dist(p,q)+\varepsilon <
\varepsilon+\varepsilon=2\varepsilon$.
\edem[Cor~\ref{inf_consequence}.]

The following lemma is useful for later constructions --- usually, it
is easier in the metric case to realize \emph{dense} subsets of
typespaces $\gaS(M)$; the lemma provides a criterion for relative
metric Galois saturation.

\begin{lemma}\label{saturation}
Suppose that we have an increasing $\ordencea$-chain of models
$(N_n:n<\omega)$ such that $N_{n+1}$ realizes a dense subset of
$\gaS(N_n)$. Then, every type in $\gaS(N_0)$ is realized in
$N_\omega:=\overline{\bigcup_{n<\omega} N_n}$.
\end{lemma}

\bdem Given $p:= \gatp(b/N_0)$ there exists $q_0\in \gaS(N_0)$
which is realized in $N_1$ (by assumption) and
$\dist(p,q_0)<\frac{1}{2(0+1)^2}=\frac12$. Let $a_0$ be a
realization of $q_0$. By corollary \ref{inf_consequence} there
exists $b_0\models p$ such that $d(b_0,a_0)<2(\frac12)=1$.
\\ \\
\indent The key idea is to build two Cauchy sequences
$(a_n)_{n<\omega}$ and $(b_n)_{n<\omega}$ such that $a_n\in N_{n+1}$,
$\gatp(b_n/N_0)=\gatp(b/N_0)$ for every $n<\omega$ and also $a_n$
and $b_n$ are closed enough, so if $c:=\lim_{n<\omega}
b_n=\lim_{n<\omega} a_n$ then by CTP
(Definition~\ref{Continuity_Types}) we have that
$\gatp(c/N_0)=\gatp(b_0/N_0)=p$. Since $c=\lim_{n<\omega}a_n$,
then $c\in N_{\omega}:=\overline{\bigcup_{n<\omega}N_n}$, and so
$p$ is realized in $N_{\omega}$.
\\ \\
\indent \underline{The construction}: Consider $n>0$. Since $N_{n+1}$ realizes a dense subset of
$\gaS(N_n)$, take $a_{n}\in N_{n+1}$ a realization of a type
$q_n\in \gaS(N_n)$ which satisfies
$\dist(\gatp(b_{n-1}/N_n),q_n)<\frac{1}{2n^2}$. By corollary
\ref{inf_consequence}, take $b_n\models \gatp(b_{n-1}/N_n)$ such
that $d(b_n,a_n)<2(\frac{1}{2n^2})=\frac{1}{n^2}$.
\\ \\
\indent We have that $(a_n)_{n<\omega}$ is a Cauchy sequence: as
$b_{n+1}\models \gatp(b_n/N_{n+1})$, there exists $g\in
Aut(\monster/N_{n+1})$ such that $g(b_n)=b_{n+1}$. Since $g$ is an
isometry and $a_n\in N_{n+1}$, then
$d(b_{n+1},a_n)=d(g(b_n),g(a_n))=d(b_n,a_n)<\frac{1}{n^2}$.
Therefore, $d(a_{n+1},a_n)\le
d(a_{n+1},b_{n+1})+d(b_{n+1},a_n)<\frac{1}{(n+1)^2}+\frac{1}{n^2}<\frac{2}{n^2}$,
so we have that $(a_n:n<\omega)$ is a Cauchy sequence.
\\ \\
\indent Therefore, there exists $c:=\lim_{n<\omega}a_n$, $c\in
N_{\omega}$ and also $c=\lim_{n<\omega}b_n$. So, we are done.
\edem[Lemma~\ref{saturation}]

\section{Smooth independence in MAECs}\label{sec:2}
Throughout this section, every model has density cardinal $\mu$ (unless
we specify a different density).
%

\begin{definition}[$\varepsilon$-splitting and $\eindep$]
Let $N\ordencea M$ and $\E>0$. We say that $\gatp(a/M)$
$\E$-\emph{splits} over $N$ iff there exist $N_1,N_2$ with  $N\ordencea N_1,N_2\ordencea
M$ and $h:N_1\cong_N N_2$ such that
$\dist(\gatp(a/N_2),h(\gatp(a/N_1))\ge \E$. We use  $a\eindep_N
M$ to denote the fact that  $\gatp(a/M)$ does
not $\varepsilon$-split over $N$,
\end{definition}



\begin{definition}\label{r_independence}
Let $N\ordencea M$. Fix $\mathcal{N}:=\langle N_i : i<\sigma \rangle$ a resolution of $N$. We say that $a$ is \emph{smooth independent} from $M$ over $N$ relative to $\C{N}$ (denoted by $a\indep^{\C{N}}_N
M$) iff for every $\E>0$ there exists $i_\E<\sigma$
such that $a\eindep_{N_{i_\E}} M$.
\end{definition}

We call {\it smooth independence} the notion of independence given above, inspired by \cite{BaSh}. In that paper, J. Baldwin and S. Shelah defined {\it smoothness} as a nice property of an abstract class of models $\C{K}$ which involves increasing chains of models, context where the existence of a kind of monster model holds.

\begin{notation}
Let $p$ be a Galois-type over $M$, $N$  a $\C{K}$-submodel of
$M$ and $\mathcal{N}$ a resolution of $N$. We denote by $p\eindep_N M$ ($p\indep_N^{\mathcal{N}} M$) iff for any
realization $a\models p$ we have that $a\eindep_N M$ ($a\indep_N^{\mathcal{N}}
M$).
\end{notation}

Next, we prove some basic properties of smooth independence.

\begin{fact}[Invariance]\label{invariance}
Let $f\in Aut(\monster)$. Then $a\indep_{M}^{\C{M}} N$ iff $f(a)\indep_{f[M]}^{f[\C{M}]} f[N]$.
\end{fact}

\begin{proposition}[Monotonicity of smooth independence]\label{monotonicity}
Let $M_0\ordencea M_1\ordencea M_2\ordencea M_3$. Fix $\mathcal{M}_k:=\langle M^k_i : i<\sigma_k \rangle$ a resolution of $M_k$ ($k=0,1$), where $\C{M}_0\subseteq \C{M}_1$. If $a\indep^{\C{M}_0}_{M_0} M_3$ then $a\indep^{\C{M}_1}_{M_1} M_2$.
\end{proposition}
\bdem
Let $\E>0$. Since $a\indep^{\C{M}_0}_{M_0} M_3$, there exists
$i_\E<\sigma_0$ such that $a\eindep_{M^0_{i_\E}} M_3$. But
$\C{M}_0\subseteq \C{M}_1$, then there exists $j_\E<\sigma_1$ such
that $M^0_{i_\E}=M^1_{j_\E}$. Therefore, for every
$M^1_{j_\E}\ordencea N_1\stackrel{h}{\cong}_{M^1_{j_\E}} N_2 \ordencea
M_3$ (in particular if $N_1,N_2\ordencea M_2$) we have that
$d(\gatp(a/N_2),\gatp(h(a)/N_2))<\E$. Then $a\eindep_{M^1_{j_\E}}
M_2$. Since this holds for every $\E>0$, then $a\indep^{\C{M}_1}_{M_1}
M_2$.
\edem[Prop.~\ref{monotonicity}]

\begin{proposition}[Monotonicity of non-$\varepsilon$-splitting]\label{monotonicity_eps}
Let $M_0\ordencea M_1\ordencea M_2\ordencea M_3$. If $a\eindep_{M_0} M_3$ then $a\eindep_{M_1} M_2$.
\end{proposition}

%
%

\begin{lemma}[Stationarity (1)]\label{Tarea1}
Suppose that $N_0\ordencea N_1 \ordencea N_2$ and $N_1$ is
universal over $N_0$. If $\gatp(a/N_1)=\gatp(b/N_1)$,
$a\eindep_{N_0} N_2$ and $b\eindep_{N_0} N_2$, then
$\dist(\gatp(a/N_2),\gatp(b/N_2)<2\varepsilon$.
\end{lemma}
\bdem Since $N_1$ is universal over $N_0$, then there exists a
$\C{K}$-embedding $g:N_2\to_{N_0} N_1$. So, $N_0\ordencea
g[N_2]\ordencea N_1$.
\\ \\
Since $N_0\ordencea g[N_2],N_2\ordencea N_2$,
$g^{-1}\upharpoonright g[N_2]: g[N_2]\cong_{N_0} N_2$ and
$a\eindep_{N_0} N_2$, then
$\dist(\gatp(g^{-1}(a)/N_2),\gatp(a/N_2) )<\varepsilon$.
\\ \\
Doing a similar argument, it is easy to prove that\linebreak
 $\dist(\gatp(g^{-1}(b)/N_2),\gatp(b/N_2) )<\varepsilon$.
\\ \\
Also, since $\gatp(a/N_1)=\gatp(b/N_1)$ and $g[N_2]\ordencea N_1$,
then $\gatp(a/g[N_2])=\gatp (b/g[N_2])$, so
$\gatp(g^{-1}(a)/N_2)=\gatp (g^{-1}(b)/N_2)$.
\\ \\
Therefore,
\begin{eqnarray*}
\dist(\gatp(a/N_2),\gatp(b/N_2)) &\le& \dist(\gatp(a/N_2),\gatp(g^{-1}(a)/N_2))\\
                                    && + \dist(\gatp(g^{-1}(a)/N_2),\gatp(g^{-1}(b)/N_2))\\
                                    && + \dist(\gatp(g^{-1}(b)/N_2),\gatp(b/N_2))\\
                                    &<& \varepsilon + 0 + \varepsilon\\
                                    &=& 2\varepsilon
\end{eqnarray*}
\edem[Lemma~\ref{Tarea1}]

\begin{proposition}[Extension of $\indep^{\C{N}}$ over universal models]\label{ExtUniv}
If $N\ordencea M \ordencea M'$, $\C{N}:=\langle N_i:i<\sigma \rangle$ is a resolution of $N$, $M$ is universal over $N$ and
$p:=\gatp(a/M)\in\gaS(M)$ is a Galois type such that $a\indep^{\C{N}}_N
M$, then there exists $b$ such that $\gatp(b/M)=\gatp(a/M)$ and
$b\indep^{\C{N}}_N M'$.
\end{proposition}
\bdem Since $M$ is universal over $N$, there exists a
$\C{K}$-embedding $h':M'\to_N M$. Extend $h'$ to an automorphism
$h\in Aut(\monster/N)$. Since $a\indep_N M$ and $h[M']\ordencea
M$, by monotonicity of $\indep^{\C{N}}$ we have that $a\indep_N h[M']$. By
invariance, we have that $h^{-1}(a)\indep^{\C{N}}_N M'$.

\begin{claim}\label{Claim}
$\gatp(a/M)=\gatp(h^{-1}(a)/M)$.
\end{claim}
\bdem Take $N_1:=h^{-1}[M]$ and $N_2:=M$. Notice that $N\ordencea
N_1,N_2\ordencea h^{-1}[M]$ and $h\upharpoonright N_1: N_1 \cong_N
N_2$. Since $a\indep^{\C{N}}_N M$, by invariance we have that
$h^{-1}(a)\indep^{\C{N}}_{N} h^{-1}[M]$. So, given
$n<\omega$ there exists $i_n<\sigma$ such that
$h^{-1}(a)\indep^{\hspace{-1mm}\frac{1}{n+1}}_{N_{i_n}} h^{-1}[M]$.
\\ \\
By monotonicity of non-$\varepsilon$-splitting
(Proposition~\ref{monotonicity_eps}), we may conclude that
$h^{-1}(a)\indep^{\hspace{-1mm}\frac{1}{n+1}}_{N} h^{-1}[M]$ for every
$n<\omega$.
\\ \\
Since $N \ordencea N_1,N_2\ordencea h^{-1}[M]$, we
have that for every $n<\omega$\linebreak
$\dist(\gatp(h^{-1}(a)/N_2),\gatp((h\circ
h^{-1})(a)/N_2))<\frac{1}{n+1}$
\\ \\
Since $N_2:=M$, we have that $\gatp(a/M)=\gatp(h^{-1}(a)/M)$. This
finishes the proof of claim \ref{Claim} \edem[Claim~\ref{Claim}]

Since $\gatp(a/M)=\gatp(h^{-1}(a)/M)$, there exists $g\in
Aut(\monster/M)$ such that $g(h^{-1}(a))=a$. Recall that
$h^{-1}(a)\indep^{\C{N}}_N M'$, so by invariance we have that
$g(h^{-1}(a))\indep^{\C{N}}_N g[M']$, i.e.: $a\indep^{\C{N}}_N g[M']$. Applying
invariance again, we have that $g^{-1}(a)\indep^{\C{N}}_N M'$. Take
$b:=g^{-1}(a)$. This now ends the proof of Proposition~\ref{ExtUniv}. \edem[Prop.~\ref{ExtUniv}]

\begin{proposition}[Stationarity (2)]\label{UnExtUniv}
If $N\ordencea M \ordencea M'$, $M$ is universal over $N$, $\C{N}:=\langle N_i:i<\sigma \rangle$ a resolution of $N$ and
$p:=\gatp(a/M)\in\gaS(M)$ is a Galois type such that $a\indep^{\C{N}}_N
M$, then there exists an unique extension $p^*\supset p$ over $M'$
which is independent (relative to $\C{N}$) from $M'$ over $N$.
\end{proposition}
\bdem By proposition \ref{ExtUniv}, there exists at least an
extension $p^{*}:=\gatp(b/M')$ of $p$ with the desired property.
\\ \\
Let $q^*:=\gatp(c/M')\supset p$ be another extension with
satisfies the desired property. So, $p^{*}\upharpoonright M =
q^*\upharpoonright M$, $b\indep^{\C{N}}_N M'$ and $c\indep^{\C{N}}_N M'$.
\\ \\
Let $\E>0$. So, there exist $i_\E^a,i_\E^b<\sigma$
such that $a\eindep_{N_{i_\E^a}} M'$ and $b\eindep_{N_{i_\E^b}}
M'$. Taking $i:=\max\{i_\E^a,i_\E^b\}$, by monotonicity of
non-$\E$-splitting we have that $a\eindep_{N_i} M'$ and
$b\eindep_{N_i} M'$.
\\ \\
Since $M$ is universal over $N_i$ (because $M$ is universal over $N$), $a\eindep_{N_i} M'$,
$b\eindep_{N_i} M'$ and $p^{*}\upharpoonright M =
q^*\upharpoonright M$, by lemma \ref{Tarea1} we have that
$\dist(p^*,q^*)<2\E$. Therefore $p^*=q^*$.
 \edem[Prop.~\ref{UnExtUniv}]

The following property of smooth independence (called {\it antireflexivity
}) is the metric version of the following property of thorn-forking in the first order setting: if $a\indep^{thorn}_B a$ then $a\in acl(B)$.

\begin{proposition}[Antireflexivity]\label{s_antireflexivity}
Let $M\ordencea N$ where $M$ is a $(\mu,\theta)$-limit model witnessed by $\C{M}:=\{M_i:i<\theta\}$. If $a\indep_{M}^{\C{M}}N$ and $a\in N$, then $a\in M$.
\end{proposition}
\bdem
Let $\E>0$ and $i_\E<\theta$ be such that $a\indep^\E_{M_{i_\E}} N$. Since $M$ is universal over $M_{i_\E}$, there exists an $\ordencea$-embedding  $f:N\to M$ which fixes $M_{i_\E}$ pointwise. Define $c:=f(a)$. Notice that $c\in M$. Setting $N_1:=N$ and $N_2:=f[N]$, notice that $M_{i_\E}\ordencea N_1 \stackrel{f}{\approx}_{M_{i_\E}} N_2\ordencea N$. Since $a\indep^\E_{M_{i_\E}} N$, then
\begin{eqnarray*}
d(\gatp(a/N_2),\gatp(c/N_2)) &=& d(\gatp(a/N_2),\gatp(f(a)/N_2))\\
 &<& \E
\end{eqnarray*}

Since $c=f(a)\in f[N]=N_2$, $c$ is the unique realization of $\gatp(c/N_2)$. Therefore, we can find $a'\models \gatp(a/N_2)$ such that $d(a',c)<\E$ (by definition of distance between types).
\\ \\
Since $a'\models \gatp(a/N_2)$, there exists $g\in Aut(\monster/N_2)$ such that $g(a)=a'$. Therefore,

\begin{eqnarray*}
d(a,c) &=& d(g(a),g(c)) \text{\ ($g$ is an isometry)}\\
&=& d(a',c) \text{\ (since $c\in N_2$ and $g\in Aut(\monster/N_2)$)}\\
&<& \E
\end{eqnarray*}

Defining  $B(a,\E):=\{b\in N : d(a,b)<\E \}$, we have that $c\in B(a,\E)\cap f[N]\subseteq B(a,\E)\cap M$, so $B(a,\E)\cap M\neq \emptyset$, hence $a\in \overline{M}=M$.
\edem[Prop. \ref{s_antireflexivity}]

\begin{proposition}[Local character of
  non-$\varepsilon$-splitting]\label{Locality}
Let $\C{K}$ be a $\mu$-$d$-stable MAEC and $\E>0$. For every $p\in
\gaS(N)$ with $N$ of density character $> \mu$ there exists $M\ordencea N$ with density character $\mu$
such that $p\eindep_M N$
\end{proposition}
\bdem Suppose that there exists some $p:=\gatp(\overline{a}/N)$ such
that $p\not\hspace{-1.3mm}\eindep_M N$ for every $M\ordencea N$
with density character $\mu$. If $\overline{a}\in N$, it is
straightforward to see that $p$ does not $\varepsilon$-split over
its domain. Then, suppose that $\overline{a}\notin N$.
\\ \\
\indent Define $\chi:=\min\{\kappa: 2^{\kappa}>\mu\}$. So,
$\chi\le \mu$ and $2^{<\chi}\le \mu$.
\\ \\
\indent We will construct a sequence of models $\langle M_\alpha,
N_{\alpha,1}, N_{\alpha,2} : \alpha<\chi \rangle$ in the following
way: First, take $M_0\ordencea N$ as any submodel of density
character $\mu$.\\ \\
\indent Suppose $\alpha:=\gamma+1$ and that $M_\gamma$ (with
density character $\mu$) has been constructed. Then $p$
$\varepsilon$-splits over $M_\gamma$. Then there exist
$M_\gamma\ordencea N_{\gamma,1},N_{\gamma,2}\ordencea N $ with
density character $\mu$ and $F_\gamma:
N_{\gamma,1}\cong_{M_\gamma} N_{\gamma,2}$ such that
$d(F_\gamma(p\rest N_{\gamma,1}),p\rest N_{\gamma,2})\ge
\varepsilon$. Take $M_{\gamma+1}\ordencea N$ a submodel of size
$\mu$ which contains $|N_{\gamma,1}|\cup |N_{\gamma,2}|$. At limit
stages $\alpha$, take
$M_\alpha:=\overline{\bigcup_{\gamma<\alpha}M_{\gamma}}$.
\\ \\
\indent Let us construct a sequence $\langle M_{\alpha}^*
:\alpha\le \chi \rangle$ of models and a tree $\langle
h_{\eta}:\eta<\alpha \rangle$ ($\alpha\le \chi$) of
$\C{K}$-embeddings such that:

\begin{enumerate}
\item $\gamma<\alpha$ implies $M_\gamma^*\ordencea M_\alpha^*$.
\item $M_{\alpha}^*:=\overline{\bigcup_{\gamma<\alpha}M_{\gamma}^*}$ if
$\alpha$ is limit.
\item $\gamma<\alpha$ and $\eta\in\;^{\alpha}2$ imply that $h_{\eta\rest \gamma}\subset
h_{\eta}$.
\item $h_\eta:M_{\alpha}\to M_{\alpha}^*$ for every $\eta\in \;^\alpha
2$.
\item If $\eta\in\; ^\gamma 2$ then $h_{\eta^\frown 0}(N_{\gamma,1})=h_{\gamma^\frown 1}(N_{\gamma,2})$
\end{enumerate}

Take $M_0^*:=M_0$ and $h_{\langle \rangle}:=id_{M_0}$.
\\ \\
\indent If $\alpha$ is limit, take
$M_{\alpha}^*:=\overline{\bigcup_{\gamma<\alpha}M_{\gamma}^*}$ and
if $\eta\in\hspace{.1mm}^\alpha 2$ define
$h_{\eta}:=\overline{\bigcup_{\gamma<\alpha}h_{\eta\rest \gamma}
}$, the unique extension of $\bigcup_{\gamma<\alpha}h_{\eta\rest
  \gamma}$ to $M_\alpha=\overline{\bigcup_{\gamma<\alpha}M_{\gamma}}$.
\\ \\
\indent If $\alpha:=\gamma+1$, let
$\eta\in\hspace{.1mm}^{\gamma}2$. Take $\overline{h_{\eta}}\supset
h_\eta$ any automorphism of the monster model $\monster$ (this is
possible because $\monster$ is homogeneous).
\\ \\
\indent Notice that $\overline{h_{\eta}}\circ
F_{\gamma}(N_{\gamma,1})=\overline{h_{\eta}}(N_{\gamma,2})$.
Define $h_{\eta^\frown 0}$ as any extension of
$\overline{h_{\eta}}\circ F_\gamma$ to $M_{\gamma+1}$ and
$h_{\eta^\frown 1}$ as $\overline{h_{\eta}}\rest M_{\gamma+1}$.
Take $M_{\gamma+1}^*\ordencea N$ as any model with density
character $\mu$ which contains $h_{\eta^\frown l}(M_{\gamma+1})$
for any $\eta\in \hspace{0.1mm}^\gamma 2$ and $l=0,1$.
\\ \\
\indent Now, for every $\eta\le\hspace{0.1mm}^\chi 2$, let $H_{\eta}$
be an automorphism of $\monster$ which extends $h_{\eta}$, 

\begin{claim}\label{LessE}
If $\eta\neq \nu \in\hspace{.1mm}^\chi 2$ then
$d(\gatp(H_\eta(\overline{a})/M_\chi^*),
\gatp(H_\nu(\overline{a})/M_\chi^*))\ge \varepsilon$.
\end{claim}
\bdem Suppose not, then $d(\gatp(H_\eta(\overline{a})/M_\chi^*),
\gatp(H_\nu(\overline{a})/M_\chi^*))< \varepsilon$. Let
$\rho:=\eta\land \nu$. Without loss of generality, suppose that $\rho^\frown 0\le \eta$
and $\rho^\frown 1\le \nu$. Let $\gamma:={\rm lg}(\rho)$. Since
$h_{\rho^\frown 0}(N_{\gamma,1})=h_{\rho^{\frown
}1}(N_{\gamma,2})\ordencea M_{\chi}^*$, then
$d(\gatp(H_\eta(\overline{a})/h_{\rho^\frown 0}(N_{\gamma,1})),
\gatp(H_\nu(\overline{a})/h_{\rho^{\frown}1}(N_{\gamma,2}))<
\varepsilon$. Also\footnote{This distance between Galois types makes sense,
  as $h_{\rho^\frown 0}(N_{\gamma,1}) = h_{\rho^\frown 1}(N_{\gamma,2})$.}
\begin{eqnarray*}
d(\gatp(H_{\nu}^{-1}\circ
H_\eta(\overline{a})/F_{\gamma}(N_{\gamma,1})),
\gatp(\overline{a}/N_{\gamma,2})) &=&\\
d(\gatp(H_\eta(\overline{a})/h_{\rho^\frown 0}(N_{\gamma,1})),
\gatp(H_\nu(\overline{a})/h_{\rho^{\frown}1}(N_{\gamma,2})) &<&
\varepsilon
\end{eqnarray*}
(as $H_{\nu}$ is an isometry, $h_{\rho^{\frown}0}=h_{\rho}\circ
F_{\gamma}$, $\rho< \nu$, $\rho^{\frown}0\le \eta$ and
$\rho^{\frown}1\le \nu$). Since $H_{\nu}^{-1}\circ
H_\eta \supset F_{\gamma}$, then
$d(F_{\gamma}(p\rest N_{\gamma,1}),p\rest N_{\gamma,2})<\E$, which
contradicts the choice of $N_{\gamma,1}$, $N_{\gamma,2}$ and
$F_{\gamma}$. This finishes the proof of claim \ref{LessE} \edem[Claim~\ref{LessE}]

We have that $dc(M_{\chi}^*)=\mu$, but claim \ref{LessE} says that
there are at least $2^{\chi}>\mu$ many types mutually at distance
at least $\E$. Therefore $dc(\gaS(M_{\chi}^*))>\mu$, which contradicts
$\mu$-$d$-stability. \edem[Prop.~\ref{Locality}]

\begin{proposition}[Existence]\label{LocalitySplitting}
Let $\C{K}$ be a $\mu$-$d$-stable MAEC. Then, for every
$\overline{a}\in \monster$ and every $N\in \C{K}$ there exists
$M\ordencea N$ with density character $\mu$ and a resolution $\C{M}:=\langle M_i : i<\omega \rangle$ of $M$ such that
$a\indep^{\C{M}}_M N$.
\end{proposition}
\bdem Let $n<\omega$. By proposition \ref{Locality}, there
exists $M_n\ordencea N$ with density character $\mu$ such that
$\overline{a}\indep_{M_n}^{\hspace{-1mm}\frac{1}{n+1}}N$. By monotonicity, without loss of generality we can
assume that $m<n<\omega$ implies $M_m\ordencea M_n$. Take
$M:=\overline{\bigcup_{n<\omega}M_n}$. Notice that
$dc(M)=\mu$. It is straightforward to see that $\overline{a}\indep^{\C{M}}_M N$.
\edem[Prop.~\ref{LocalitySplitting}]

\begin{lemma}[Continuity of independence]\label{cont_indep}
Let $(b_n)_{n<\omega}$ be a convergent sequence and $b:=\lim_{n<\omega} b_n$. If $b_n\indep^{\C{N}}_N M$ for every $n<\omega$, then $b\indep^{\C{N}}_N M$.
\end{lemma}
\bdem
Since $b_n\indep^{\C{N}}_N M$ ($n<\omega$), for every $\varepsilon>0$ there exists $i_{n,\E}<\sigma$ such that for every $N_{i_{n,\E}}\ordencea N^1\stackrel{h}{\cong}_{N^n_{\E}}N^2 \ordencea M$, therefore we have that $d(\gatp(b_n/N^2),\gatp(h(b_n)/N^2))<\varepsilon/3$.
\\ \\
Let $K<\omega$ be such that for every $n\ge K$  we have that $d(b_n,b)<\varepsilon/3$. Therefore, $d(\gatp(b_n/N^2),\gatp(b/N^2))<\varepsilon/3$ for every $n\ge K$.
\\ \\
Since $h$ is an isometry, we have that $(h(b_n))\to h(b)$ and also for every $n\ge K$ we have that $d(h(b_n),h(b))<\varepsilon/3$ (and therefore\linebreak
$d(\gatp(h(b_n)/N^2),\gatp(h(b)/N^2))<\varepsilon/3$).
\\ \\
\noindent Therefore, for any $n\ge K$ we have that
\begin{eqnarray*}
d(\gatp(h(b)/N^2),\gatp(b/N^2)) &\le& d(\gatp(h(b)/N^2),\gatp(h(b_n)/N^2))+\\
&\ & d(\gatp(h(b_n)/N^2),\gatp(b_n/N^2))+\\
&\ & d(\gatp(b_n/N^2),\gatp(b/N^2))\\
&<& \varepsilon/3+\varepsilon/3+\varepsilon/3=\varepsilon.
\end{eqnarray*}
Therefore, $b\indep^\varepsilon_{N_{i_{n,\E}}}M$ and so, $b\indep^{\C{N}}_N M$.
\edem[Lemma~\ref{cont_indep}]

\begin{proposition}[stationarity (3)]\label{Stationarity}
Let $M_0\ordencea M\ordencea N$ be such that $M$ is a $(\mu,\sigma)$-limit model over $M_0$, where $\C{M}:=\{M_i:i<\sigma\}$ witnesses that $M$ is $(\mu,\sigma)$-limit over $M_0$. If $a,b\indep^{\C{M}}_M N$ and $\gatp(a/M)=\gatp(b/M)$, then $\gatp(a/N)=\gatp(b/N)$.
\end{proposition}
\bdem
Let $\E>0$. Since $a,b\indep^{\C{M}}_M N$, there exists $i<\sigma$ such that $a,b\eindep_{M_i}N$ (by definition and monotonicity of non-$\E$-splitting). Since $M_{i+1}$ is universal over $M_i$ and $M_i\ordencea N$, there exists and $\ordencea$-embedding $f:N\to_{M_i} M_{i+1}$.  Also, since $M_i\ordencea f[N]\stackrel{f^{-1}}{\cong_{M_i}} N \ordencea N$ and $a\eindep_{M_i} N$, therefore $d(\gatp(a/N),\gatp(f^{-1}(a)/N))<\E$.
\\ \\
\indent Doing a similar argument, we have that $d(\gatp(b/N),\gatp(f^{-1}(b)/N))<\E$.
\\ \\
\indent On the other hand, we have that $\gatp(a/f[N])=\gatp(b/f[N])$ (since $\gatp(a/M)=\gatp(b/M)$ and $f[N]\ordencea M_{i+1}\ordencea M$), therefore we have that $\gatp(f^{-1}(a)/N)=\gatp(f^{-1}(a)/N)$.
\\ \\
\indent Hence
\begin{eqnarray*}
d(\gatp(a/N),\gatp(b/N)) &\le& \gatp(a/N)+\gatp(f^{-1}(a)/N)\\
&+& d(\gatp(f^{-1}(a)/N),\gatp(f^{-1}(b)/N))\\
&+& d(\gatp(f^{-1}(b)/N),\gatp(b/N))\\
&<& \E + 0 + \E\\
&=& 2\E
\end{eqnarray*}

Therefore, $\gatp(a/N)=\gatp(b/N)$.
\edem[Prop.~\ref{Stationarity}]

\begin{proposition}[Transitivity]\label{transitivity}
Let $M_0\ordencea M_1\ordencea M_2$ be such that $M_1$ and $M_0$ are $(\mu,\sigma)$-limit over some $M'\ordencea M_0\ordencea M_1$, where $\C{M}_i$ witnesses that $M_i$ is $(\mu,\sigma)$-limit over $M'$ and $\C{M}_0\subset \C{M}_1$. Then $a\indep^{\C{M}_0}_{M_0} M_2$ iff $a\indep^{\C{M}_0}_{M_0} M_1$ and $a\indep^{\C{M}_1}_{M_1} M_2$.
\end{proposition}
\bdem
$(\Rightarrow)$ By monotonicity.
\\ \\
\indent  $(\Leftarrow)$ Suppose $a\indep^{\C{M}_0}_{M_0} M_1$ and $a\indep^{\C{M}_1}_{M_1} M_2$.
Notice that $M_1$ is universal over $M_0$. Therefore, by extension property (proposition \ref{ExtUniv}), there exists $b\models \gatp(a/M_1)$ such that $b\indep^{\C{M}_0}_{M_0} M_2$. By monotonicity, we have that $b\indep^{\C{M}_1}_{M_1} M_2$. Since $a,b\indep^{\C{M}_1}_{M_1} M_2$, $\gatp(a/M_1)=\gatp(b/M_1)$ and in particular $M_1$ is a limit model over $M_0$, then by stationarity (proposition \ref{Stationarity}) we have that $\gatp(a/M_2)=\gatp(b/M_2)$. Since $b\indep^{\C{M}_0}_{M_0} M_2$, then $a\indep^{\C{M}_0}_{M_0} M_2$.
\edem[Prop.~\ref{transitivity}]

\section[Domination, orthogonality and parallelism
$\dots$]{Domination, orthogonality and parallelism under uniqueness of
  limit models}\label{chapter:3}

The study of  Zilber's trichotomy for strongly minimal sets in
understanding the classification -up to bi-interpretability- of
uncountably categorical strongly minimal theories is an important
step toward geometric stability theory (although restricted to
$\omega_1$-categoricity).
The non-finite
axiomatizability of totally categorical theories -works of Cherlin,
Harrington, Lachlan and Zilber- is the main initial step toward
geometric stability theory. Buechler used generalizations of this
machinery outside of totally-categorical and $\omega_1$-categorical
settings and obtained a proof of his famous dichotomy theorem on the
collection $D$ of realizations of $stp(a/A)$ for $a$ any realization
of a weakly minimal type, which says that either $D$ is locally
modular or $p$ has Morley rank $1$  \cite{Bu2}.

 In Superstable First Order theories, the development of
 Geometric Stability Theory (see \cite{Pi,Bu}) includes
 generalizations of results studied in the categorical settings. In
 his doctoral thesis, E. Hrushovski extended this work to Stable First
 Order Theories \cite{Hr}. Also, this study has been extended to Rosy
 Theories by A. Onshuus and A. Usvyatsov (see \cite{OnUs1}). In
 abstract settings, S. Shelah provided  some extensions of these
 results in {\it AEC Good Frames} (see \cite{Sh705}), which
 corresponds to a setting that J. Baldwin calls {\it intermediate
   stability theory} because it does not really consider more refined
 techniques of geometric stability theory, e.g. group configurations
 and Hrushovski's analysis.

This chapter is devoted to the study of some basic geometrical notions
of classical stability theory: {\it domination, orthogonality and
  parallelism}. These notions correspond in the MAEC setting to the
well-known notions going by the same names in stable first order
theories. We will study some of their properties in MAEC settings
exhibiting behavior akin to (variants of) superstability, and will extend
results due to Baldwin
(\cite{Ba}) and Shelah \cite{Sh705}.


\begin{assumption}\label{Superstability}
Throughout this section, we assume AP, JEP, CTP, existence of
arbitrarily large enough models and  the following assumptions (we
sometimes abusively call them {\it ``superstability''} - but we do not
attempt to define that notion at this stage):
For every $a$ and every increasing and continuous $\ordencea$-chain of models
$\langle M_i : i<\sigma \rangle$ and $M _j$ a resolution of $M_j$ ($j<\sigma$):
\begin{enumerate}
\item (Continuity) If $p\rest M_i \indep^{\C{M} _0}_{M_0} M_i$ for all $i<\sigma$, then
$p\indep^{\C{M} _0}_{M_0} \overline{\bigcup_{i<\sigma} M_i}$.
\item (Locality) if $cf(\sigma)>\omega$, there exists $j<\sigma$ such that
$a\indep^{\C{M} _j}_{M_j} \bigcup_{i<\sigma} M_i$.
\item ($\E$-simplicity) if $cf(\sigma)=\omega$, there exists $j<\sigma$ such that
$a\indep^{\E}_{M_j} \overline{\bigcup_{i<\sigma} M_i}$.
\end{enumerate}
\end{assumption}

\begin{remark}\label{Simplicity}
Under these assumptions plus categoricity, we proved in~\cite{ViZa}
the uniqueness -up to
isomorphisms- of limit models over a fixed base: If
$M_i$ is a
$(\mu,\theta_i)$-$d$-limit over $M$ ($i\in \{1,2\}$) such that
$dc(M_1)=dc(M_2)$, then $M_1\approx_M M_2$. It is straightforward to
see that assumptions \ref{Superstability} 2. and 3. imply
$a\indep^{\C{M}}_M M$ for every $M$ and every resolution $\C{M}$ of
$M$.
\end{remark}

\begin{assumption}[Uniqueness of limit models]\label{Uniqueness}
If $M$ and $N$ are limit models of the same density character $\mu$ over the same model $M_0$, then there exists an isomorphism $f:M\to N$ fixing $M_0$ pointwise.
\end{assumption}

\subsection{Domination in MAEC}\label{section:domination}
In this section, we define a natural adaptation of the notion of domination in the setting of superstable MAECs that exhibit the superstability-like assumption \ref{Superstability}. We base the development of this section on \cite{Ba} but we use s-independence instead of intersections as Baldwin does.
\\ \\
\indent According to S. Buechler (\cite{Bu}), the motivating question which takes us to the notion of {\it domination} is whether nonorthogonal (first order syntactical) types $p$ and $q$ have bases relative to a model $M$ (i.e., maximal Morley sequences of $p$ and $q$ respectively over the domain of the respective types contained in $M$) with the same cardinality. In such context, domination is a kind of opposite notion to orthogonality. In first order, we say that a (possibly infinite) set $B$ dominates another (possibly infinite) set $A$ over $C$ if and only if for any set $D$, if $B\indep_C D$ then $A\indep_C D$. But in our setting, we cannot define independence on sets because, in general, Galois types are defined on models. Because of that, we have to adapt this notion to our general context.  

\begin{notation}
$(M,\C{M},N,a)$ means that $M\ordencea N$, $M$ is a limit model witnessed by $\C{M}$ and $a\in N\setminus M$.
\end{notation}

\begin{definition}
We say that $(M,\C{M},N,a)\prec_{nf} (M',\C{M}',N',a)$ if and only if $M'$ is a limit model over $M$, $\C{M}\subset \C{M}'$ and $\C{M}$ corresponds to an initial segment of $\C{M}'$, $N\ordencea N'$ and $a\indep^{\C{M}}_M M'$.
\end{definition}

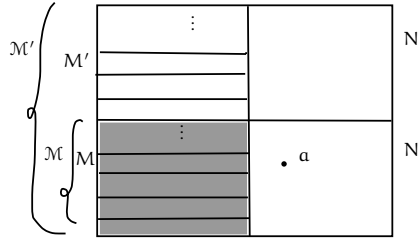
\begin{center}
\scalebox{0.6} 
{
\begin{pspicture}(0,-2.6)(9.182813,2.62)
\definecolor{color63b}{rgb}{0.6,0.6,0.6}
\psframe[linewidth=0.04,dimen=outer](8.320937,2.56)(1.7409375,-2.6)
\psline[linewidth=0.04cm](5.1009374,2.52)(5.1409373,-2.58)
\psline[linewidth=0.04cm](1.7809376,0.0)(8.280937,-0.02)
\psframe[linewidth=0.04,linecolor=color63b,dimen=outer,fillstyle=solid,fillcolor=color63b](5.0809374,-0.06)(1.8209375,-2.54)
\psline[linewidth=0.04cm,fillcolor=color63b](1.7809376,-2.2)(5.1409373,-2.18)
\psline[linewidth=0.04cm,fillcolor=color63b](1.7809376,-1.72)(5.1009374,-1.72)
\psline[linewidth=0.04cm,fillcolor=color63b](1.8009375,-1.18)(5.1009374,-1.18)
\psline[linewidth=0.04cm,fillcolor=color63b](1.7609375,-0.76)(5.1409373,-0.74)
\psline[linewidth=0.04cm,fillcolor=color63b](1.7809376,0.46)(5.0809374,0.46)
\psline[linewidth=0.04cm,fillcolor=color63b](1.7209375,1.0)(5.0409374,1.02)
\psline[linewidth=0.04cm,fillcolor=color63b](1.8209375,1.5)(5.0609374,1.46)
\psline[linewidth=0.04cm,fillcolor=color63b](5.0609374,1.46)(5.0609374,1.5)
\rput(1.5,-0.85){$M$}
\rput(1.3,1.33){$M'$}
\rput(8.792344,1.81){$N'$}
\rput(8.732344,-0.65){$N$}
\psdots[dotsize=0.12](5.9009376,-0.98)
\rput(6.342344,-0.81){$a$}
\rput(3.8823438,2.25){$\vdots$}
\rput(3.6423438,-0.17){$\vdots$}
\pscustom[linewidth=0.04]
{
\newpath
\moveto(1.4409375,-0.04)
\lineto(1.4009376,-0.01)
\curveto(1.3809375,0.0050)(1.3259375,-0.01)(1.2909375,-0.04)
\curveto(1.2559375,-0.07)(1.2209375,-0.235)(1.2209375,-0.37)
\curveto(1.2209375,-0.505)(1.2209375,-0.745)(1.2209375,-0.85)
\curveto(1.2209375,-0.955)(1.2209375,-1.125)(1.2209375,-1.19)
\curveto(1.2209375,-1.255)(1.2109375,-1.39)(1.2009375,-1.46)
\curveto(1.1909375,-1.53)(1.1309375,-1.635)(1.0809375,-1.67)
\curveto(1.0309376,-1.705)(0.9709375,-1.69)(0.9609375,-1.64)
\curveto(0.9509375,-1.59)(0.9659375,-1.535)(0.9909375,-1.53)
\curveto(1.0159374,-1.525)(1.0609375,-1.54)(1.0809375,-1.56)
\curveto(1.1009375,-1.58)(1.1159375,-1.665)(1.1109375,-1.73)
\curveto(1.1059375,-1.795)(1.0959375,-1.925)(1.0909375,-1.99)
\curveto(1.0859375,-2.055)(1.0959375,-2.145)(1.1109375,-2.17)
\curveto(1.1259375,-2.195)(1.1659375,-2.235)(1.1909375,-2.25)
\curveto(1.2159375,-2.265)(1.2709374,-2.285)(1.3009375,-2.29)
}
\pscustom[linewidth=0.04]
{
\newpath
\moveto(0.8209375,2.6)
\lineto(0.7409375,2.51)
\curveto(0.7009375,2.465)(0.6209375,2.26)(0.5809375,2.1)
\curveto(0.5409375,1.94)(0.5009375,1.67)(0.5009375,1.56)
\curveto(0.5009375,1.45)(0.5009375,1.21)(0.5009375,1.08)
\curveto(0.5009375,0.95)(0.4609375,0.64)(0.4209375,0.46)
\curveto(0.3809375,0.28)(0.3159375,0.075)(0.2909375,0.05)
\curveto(0.2659375,0.025)(0.2359375,0.045)(0.2309375,0.09)
\curveto(0.2259375,0.135)(0.2259375,0.2)(0.2309375,0.22)
\curveto(0.2359375,0.24)(0.2759375,0.26)(0.3109375,0.26)
\curveto(0.3459375,0.26)(0.3859375,0.195)(0.3909375,0.13)
\curveto(0.3959375,0.065)(0.3959375,-0.15)(0.3909375,-0.3)
\curveto(0.3859375,-0.45)(0.3809375,-0.68)(0.3809375,-0.76)
\curveto(0.3809375,-0.84)(0.3659375,-1.07)(0.3509375,-1.22)
\curveto(0.3359375,-1.37)(0.3209375,-1.59)(0.3209375,-1.66)
\curveto(0.3209375,-1.73)(0.3209375,-1.885)(0.3209375,-1.97)
\curveto(0.3209375,-2.055)(0.3359375,-2.16)(0.3509375,-2.18)
\curveto(0.3659375,-2.2)(0.4059375,-2.25)(0.4309375,-2.28)
\curveto(0.4559375,-2.31)(0.5059375,-2.365)(0.5309375,-2.39)
\curveto(0.5559375,-2.415)(0.6159375,-2.455)(0.6509375,-2.47)
\curveto(0.6859375,-2.485)(0.7459375,-2.505)(0.7709375,-2.51)
\curveto(0.7959375,-2.515)(0.8659375,-2.52)(0.9109375,-2.52)
\curveto(0.9559375,-2.52)(1.0159374,-2.52)(1.0609375,-2.52)
}
\rput(0.1,1.59){$\C{M'}$}
\rput(0.8,-0.75){$\C{M}$}
\end{pspicture}
}
\end{center}

\begin{definition}
We say that a set $A$ is smooth independent from $M$ over $N$ relative to a resolution $\C{N}$ of $N$ -denoted by $A\indep^{\C{N}}_N M$ - if and only if $b\indep^{\C{N}}_N M$ for every finite tuple $b\in A$.
\end{definition}

\begin{definition}
Given $(M,\C{M},N,a)$, we say that $a$ {\it dominates} $N$ over $M$ relative to $\C{M}$ (denoted by $a\dominates^{\C{M}}_M N$) iff for every $(M',\C{M}',N',a)\succ_{nf} (M,\C{M},N,a)$ we have that $N\indep^{\C{M}}_M M'$ (i.e., for every $b\in N$ $b\indep^{\C{M}}_M M')$.
\end{definition}

Remember that in first order, $B$ dominates $A$ over $C$ if and only if for any set $D$, if $B\indep_C D$ then $A\indep_C D$. Because in our general context Galois types are defined on models instead of sets, we have to adapt this notion to our setting. Notice that $(M',\C{M}',N',a)\succ_{nf} (M,\C{M},N,a)$ implies $a\indep_{M}^{\C{M}} M'$, so $a\dominates_M^{\C{M} } N$ means that $a\indep_M^{\C{M}} M'$ implies $N\indep_M^{\C{M}} M'$, agreeing with the first order notion of domination.  

%

The following proposition says that domination over a model $M_\alpha$ implies domination over a  $\C{K}$-superstructure $M\ordenceac M_\alpha$ if there is some independence from $M$ over $M_\alpha$ (i.e., the information given over $M$ is the same over $M_\alpha$).

\begin{proposition}\label{Domination_up}
Let $(M,\C{M},N,a)$ (where $\C{M}:=\{M_i:i<\theta\}$ witnesses that $M$ is a limit model) and $\C{M}_\alpha\subset \C{M}$ be a resolution of $M_\alpha$ ($\alpha<\theta$) such that $a\indep_{M_\alpha}^{\C{M}_\alpha} M$. If $a\dominates_{M_\alpha}^{\C{M}_\alpha} N$ then $a\dominates_{M}^{\C{M}} N$.   
\end{proposition}
\bdem
Let $(M',\C{M}',N',a)\succ_{nf} (M,\C{M},N,a)$. Therefore, $a\indep_{M}^{\C{M}}M'$. By hypothesis $a\indep_{M_\alpha}^{\C{M}_\alpha} M$, hence $a\indep_{M_{\alpha}}^{\C{M}_\alpha} M'$ (by transitivity, proposition \ref{transitivity}). So, $(M',\C{M}',N',a)\succ_{nf} (M_\alpha,\C{M}_\alpha,N,a)$. Since $a\dominates_{\C{M_\alpha}}^{\C{M}_\alpha} N$, then $N\indep_{M_\alpha}^{\C{M}_\alpha} M'$. By monotonicity (proposition \ref{monotonicity}), $N\indep_{M}^{\C{M}}M'$, therefore $a\dominates_M^{\C{M}} N$.
\edem[Prop. \ref{Domination_up}]

The following proposition is a kind of reciprocal of proposition~\ref{Domination_up}. This says that under some independence from $M$ over $M_\alpha$, domination over $M$ implies domination over $M_\alpha$.  

\begin{proposition}\label{Domination_down}
Let $(M,\C{M},N,a)$ (where $\C{M}:=\{M_i:i<\theta\}$ witnesses that $M$ is a $(\mu,\sigma)$-limit model) and $\C{M}_\alpha\subset \C{M}$ be a resolution of $M_\alpha$ ($\alpha<\theta$) such that $N\indep_{M_\alpha}^{\C{M}_\alpha} M$. If $a\dominates_{M}^{\C{M}} N$ then $a\dominates_{M_\alpha}^{\C{M}_\alpha} N$. 
\end{proposition}
\bdem
Let $(M',\C{M}',N',a)\succ_{nf} (M_\alpha,\C{M}_\alpha,N,a)$. Let $M\cup M'\subset \hat{M}\ordencea \monster$  (by downward L\"owenheim-Skolem axiom) and $M^*\ordenceac \hat{M}$ be a limit over $\hat{M}$ -and so $M^*$ is a limit model over $M$, where $\C{M}^{**}$ is a witness of that-. Let $N^*\ordenceac N$ be such that $N^*\ordenceac M^*$. and $\C{M^*}:=\C{M}^\frown \C{M}^{**}$
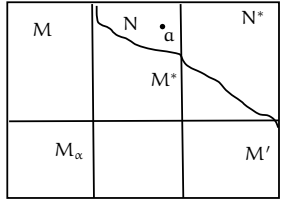
\begin{center}
\scalebox{0.6} 
{
\begin{pspicture}(0,-2.19)(6.06,2.19)
\psframe[linewidth=0.04,dimen=outer](6.06,2.19)(0.0,-2.19)
\psline[linewidth=0.04cm](1.9,2.13)(1.94,-2.17)
\psline[linewidth=0.04cm](0.06,-0.47)(6.04,-0.45)
\rput(1.39,-1.125){$M_\alpha$}
\rput(0.8,1.555){$M$}
\psline[linewidth=0.04cm](3.86,2.13)(3.92,-2.17)
\rput(2.75,1.695){$N$}
\psdots[dotsize=0.12](3.46,1.63)
\rput(3.6,1.4){$a$}
\rput(5.57,-1.165){$M'$}
\rput(5.46,1.835){$N^*$}
\pscustom[linewidth=0.04]
{
\newpath
\moveto(2.0,2.09)
\lineto(2.0,1.94)
\curveto(2.0,1.865)(2.015,1.76)(2.03,1.73)
\curveto(2.045,1.7)(2.12,1.655)(2.18,1.64)
\curveto(2.24,1.625)(2.32,1.575)(2.34,1.54)
\curveto(2.36,1.505)(2.425,1.45)(2.47,1.43)
\curveto(2.515,1.41)(2.59,1.385)(2.62,1.38)
\curveto(2.65,1.375)(2.75,1.315)(2.82,1.26)
\curveto(2.89,1.205)(3.175,1.125)(3.39,1.1)
\curveto(3.605,1.075)(3.835,1.03)(3.85,1.01)
\curveto(3.865,0.99)(3.89,0.94)(3.9,0.91)
\curveto(3.91,0.88)(3.95,0.82)(3.98,0.79)
\curveto(4.01,0.76)(4.09,0.71)(4.14,0.69)
\curveto(4.19,0.67)(4.28,0.63)(4.32,0.61)
\curveto(4.36,0.59)(4.42,0.555)(4.44,0.54)
\curveto(4.46,0.525)(4.53,0.49)(4.58,0.47)
\curveto(4.63,0.45)(4.7,0.41)(4.72,0.39)
\curveto(4.74,0.37)(4.795,0.305)(4.83,0.26)
\curveto(4.865,0.215)(4.98,0.12)(5.06,0.07)
\curveto(5.14,0.02)(5.235,-0.05)(5.25,-0.07)
\curveto(5.265,-0.09)(5.32,-0.135)(5.36,-0.16)
\curveto(5.4,-0.185)(5.475,-0.235)(5.51,-0.26)
\curveto(5.545,-0.285)(5.605,-0.315)(5.63,-0.32)
\curveto(5.655,-0.325)(5.705,-0.33)(5.73,-0.33)
\curveto(5.755,-0.33)(5.8,-0.34)(5.82,-0.35)
\curveto(5.84,-0.36)(5.88,-0.385)(5.9,-0.4)
\curveto(5.92,-0.415)(5.955,-0.455)(5.97,-0.48)
\curveto(5.985,-0.505)(6.005,-0.55)(6.02,-0.61)
}
\rput(3.5,0.455){$M^*$}
\end{pspicture} 
}
\end{center}
Since $a\indep_{M_\alpha}^{\C{M}_\alpha} M'$ (by definition of $\prec_{nf}$) and $M'$ is universal over $M_\alpha$, by the extension property of smooth independence (proposition~\ref{ExtUniv}), there exists $a'\models \gatp(a/M')$ such that $a'\indep_{M_\alpha}^{\C{M}_\alpha} M^*$. Without less of generality, suppose $a\indep_{M_\alpha}^{\C{M}_\alpha} M^*$. Notice that $(M^*,\C{M}^*,N^*,a)\succ_{nf} (M,\C{M},N,a)$. Since $a\dominates_{M}^{\C{M}} N$, then $N\indep_M M^*$. By hypothesis, $N\indep_{M_\alpha}^{\C{M}_{\alpha}} M$, so by transitivity (proposition~\ref{transitivity}, since $M$ and $M_\alpha$ are limit models over $M_0$) $N\indep_{M_\alpha} M^*$, and by monotonicity (proposition~\ref{monotonicity}) $N\indep_{M_\alpha}^{\C{M}_{\alpha}} M'$ (since $M_\alpha\ordencea M'\ordencea M^*$). So, we have that $a\dominates^{\C{M}_\alpha}_{M_\alpha} N$.
\edem[Prop. \ref{Domination_down}]

The following proposition says that given any tuple $(M,\C{M},N,a)$, we can find some extensions $N'\ordenceac N$ and $M'\ordenceac M$ such that $N'\ordenceac M'$, and $a$ dominates $N'$ over $M'$.

\begin{proposition}\label{Extension_domination}
Given $(M,\C{M},N,a)$ there exists 
$(M',\C{M}',N',a)\succ_{nf}(M,\C{M},N,a)$ such that $a\dominates^{\C{M}'}_{M'}N'$.
\end{proposition}
\bdem
Suppose not. This allows us to construct an $\prec_{nf}$-increasing and continuous sequence of tuples $\langle (M^\alpha, \C{M}^\alpha,N^\alpha,a) : \alpha<\mu^+ \rangle$ such that $(M^0,\C{M}^0,N^0,a):=(M,\C{M},N,a)$ and 
$(M^{\alpha+1},\C{M}^{\alpha+1},N^{\alpha+1},a)$ witnesses that $(M^\alpha,\C{M}^{\alpha},N^\alpha,a)$ does not satisfy that $a\dominates^{\C{M}^\alpha}_{M^\alpha} N^{\alpha}$. Therefore, there exists $b\in N^\alpha$ such that $a\indep^{\C{M}^\alpha}_{M^\alpha} M^{\alpha+1}$ but $b\not\hspace{-2mm}\indep^{\C{M}^{\alpha}}_{M^\alpha} M^{\alpha+1}$. By assumption \ref{Superstability}, given any $c$ there exists $\alpha_c<\mu^+$ such that $c\indep^{\C{M}^{\alpha_c}}_{M^{\alpha_c}} \bigcup_{\alpha<\mu^+} M^{\alpha}$.
\\ \\
\indent Consider $\gamma_0<\mu^+$. Since $N^{\gamma_0}$ has density character $\mu$, there exists $B_{\gamma_0}$ a dense subset of $ N^{\gamma_0}$ of cardinality $\mu$. Defining $f_0:B_{\gamma_0}\to \mu^+$ as $f(c):=\alpha_c$, we have that there exists $\gamma_0'<\mu^+$ such that $f(c):=\alpha_c<\gamma_0'$ for every $c\in B_{\gamma_0}$. Define $\gamma_1:=\max\{\gamma_0,\gamma_0'\}+1$.
\\ \\
\indent In the same way we define $B_{\gamma_n}$ and $\gamma_{n}$ for every $n<\omega$. Notice that $(\gamma_n:n<\omega)$ is an increasing sequence of ordinals $<\mu^+$.
\\ \\
\indent Define $\gamma:=\sup\{\gamma_n : n<\omega\}$. Notice that $\gamma<\mu^+$.
\\ \\
\indent Let $b\in N^\gamma$ be such that $b\not\hspace{-2mm}\indep^{\C{M}^\gamma}_{M^{\gamma}} M^{\gamma+1}$. Since $N^{\gamma}:=\overline{\bigcup_{\alpha<\gamma} N^\alpha}$, there exists a sequence $(b_n)\in \bigcup_{\alpha<\gamma} N^\alpha$ such that $(b_n)\to b$. By proposition \ref{cont_indep} (continuity of $\indep$), there exists $k<\omega$ such that $b_k\not\hspace{-2mm}\indep^{\C{M}^\gamma}_{M^{\gamma}} M^{\gamma+1}$. Since $b_k\in \bigcup_{\alpha<\gamma} N^\alpha$, there exists $\beta<\gamma$ such that $b_k\in N^\beta$. Since $\beta<\gamma:=\sup\{\gamma_n:n<\omega\}$, there exists $m<\omega$ such that $\beta<\gamma_m$, so $b_k\in N^{\gamma_m}$. Since by construction we have that $\overline{B_{\gamma_m}}=N^{\gamma_m}$, there exists a sequence $(c_n)\in B_{\gamma_m}$ such that $(c_n)\to b_k$. By proposition \ref{cont_indep} again, there exists $l<\omega$ such that $c:=c_l\not\hspace{-2mm}\indep^{\C{M}^\gamma}_{M^{\gamma}} M^{\gamma+1}$. By construction, $\alpha_c<\gamma_{m+1}<\gamma<\gamma+1<\mu^+$, then by proposition \ref{monotonicity} (monotonicity of $\indep$) we have that $c\not\hspace{-2mm}\indep^{\C{M}^{\alpha_c}}_{M^{\alpha_c}} \bigcup_{\alpha<\mu^+}M^{\alpha}$ (contradiction). Therefore, the proposition is true.
\edem[Prop. \ref{Extension_domination}]

The following proposition says that under the conclusions of the previous proposition, we can find an extension $N^*$ of $N'$ such that $a$ dominates $N^*$ over $M$.  

\begin{proposition}\label{EquiDomination}
Suppose $(M,\C{M},N,a)\prec_{nf} (M',\C{M}',N',a)$, where $M$  is a $(\mu,\sigma_1)$-limit model witnessed by $\C{M}:=\langle M_i:i<\sigma_1\rangle$, $M'$ is a $(\mu,\sigma_2)$-limit model over $M$ witnessed by $\C{M}''$ and $\C{M}':=\C{M}^\frown \C{M}'$, $a\indep^{\C{M}_\alpha}_{M_\alpha} M$ for some limit $\alpha<\sigma$ and $a\dominates^{\C{M}'}_{M'} N'$. Then, there exist $N^*$ and a resolution $\C{M}^*$ which witnesses that $M$ is a limit model over $M_0$ such that $a\dominates^{\C{M}^*}_M N^*$.
\end{proposition}
\bdem
Let $p:=\gatp(a/M)$ and $p':=\gatp(a/M')$. Since $a\indep^{\C{M}_\alpha}_{M_\alpha} M$ (by\linebreak hypothesis), $a\indep^{\C{M}}_{M}M'$ (by definition of $\prec_{nf}$) and $M$, $M_\alpha$ are limit models over $M_0$, by transitivity (proposition \ref{transitivity}) we have  $a\indep^{\C{M}_\alpha}_{M_\alpha}M'$.

Notice that since $M$ and $M'$ are limit over $M_\alpha$ witnessed by
$\C{M}$ and $\C{M}'$ respectively such and $\C{M}\subset \C{M}'$, then
$M$ and $M'$ are limit over $M_1\in \C{M}$. By
assumption~\ref{Uniqueness} (uniqueness of limit models), there exists
\mbox{$f:M'\stackrel{\approx}{\to}_{M_{\alpha+1}} M$}. Since
$a\indep^{\C{M}_\alpha}_{M_\alpha} M'$, we have that
$f(a)\indep^{\C{M}_\alpha}_{M_\alpha} M$ (by invariance, proposition
\ref{invariance}). 
Notice that $M_{\alpha+1}$ is universal over $M_\alpha$. Then, as $\gatp(a/M_{\alpha+1})=\gatp(f(a)/M_{\alpha+1})$  and $a,f(a)\indep^{\C{M}_\alpha}_{M_\alpha} M$,  by stationarity (proposition \ref{UnExtUniv}) we may say $\gatp(a/M)=\gatp(f(a)/M)$.

Consider $g\in Aut(\monster/M)$ such that $(g\circ f )(a)=a$. Notice that

$$
(g\circ f)(M',\C{M}',N',a) = (M,(g\circ f)[\C{M}'],(g\circ f)[N'],a)
$$
witnesses that $a\dominates^{\C{M}^*}_{M} N^*$, where $N^*:=(g\circ f)[N']$ and $\C{M}^*:=(g\circ f)[\C{M}']=f[\C{M}']$. Notice that $\C{M}^*$ is also a resolution which witnesses that $M$ is a limit model over $M_0$ (remember that in particular $f$ fixes $M_0$ pointwise).

\begin{center}
\scalebox{0.6} 
{
\begin{pspicture}(0,-2.16)(7.48,2.14)
\psframe[linewidth=0.04,dimen=outer](5.56,2.14)(0.58,-2.12)
\psline[linewidth=0.04cm](3.06,2.06)(3.06,-2.14)
\psline[linewidth=0.04cm](0.62,-0.1)(5.58,-0.1)
\psdots[dotsize=0.12](3.84,1.04)
\rput(4.2,1.125){$a$}
\rput(5.98,1.885){$N'$}
\rput(5.2,-1.8){$N$}
\psline[linewidth=0.04cm,linestyle=dashed,dash=0.16cm 0.16cm,doubleline=true,doublesep=0.12](0.6,-1.72)(3.0,-1.72)
\psline[linewidth=0.04cm,linestyle=dashed,dash=0.16cm 0.16cm,doubleline=true,doublesep=0.12](0.6,-1.3)(3.02,-1.3)
\psline[linewidth=0.04cm,linestyle=dashed,dash=0.16cm 0.16cm,doubleline=true,doublesep=0.12](0.58,-0.86)(3.0,-0.84)
\psline[linewidth=0.04cm,linestyle=dashed,dash=0.16cm 0.16cm,doubleline=true,doublesep=0.12](0.6,-0.36)(2.94,-0.4)
\psline[linewidth=0.04cm,linestyle=dashed,dash=0.16cm 0.16cm,doubleline=true,doublesep=0.12](0.58,0.52)(3.0,0.48)
\psline[linewidth=0.04cm,linestyle=dashed,dash=0.16cm 0.16cm,doubleline=true,doublesep=0.12](0.56,0.2)(2.98,0.16)
\psline[linewidth=0.04cm,linestyle=dashed,dash=0.16cm 0.16cm,doubleline=true,doublesep=0.12](0.6,0.9)(3.08,0.84)
\psline[linewidth=0.04cm,linestyle=dashed,dash=0.16cm 0.16cm,doubleline=true,doublesep=0.12](0.62,1.28)(3.02,1.26)
\psline[linewidth=0.04cm,linestyle=dashed,dash=0.16cm 0.16cm,doubleline=true,doublesep=0.12](0.62,1.66)(3.02,1.66)
\psline[linewidth=0.04cm,linestyle=dashed,dash=0.16cm 0.16cm,doubleline=true,doublesep=0.12](0.64,2.0)(3.04,1.98)
\rput(0.2,1.705){$M'$}
\rput(0.15,-0.535){$M$}
\pscustom[linewidth=0.04,linestyle=dashed,dash=0.16cm 0.16cm]
{
\newpath
\moveto(3.06,-0.04)
\lineto(3.12,0.2)
\curveto(3.15,0.32)(3.185,0.48)(3.19,0.52)
\curveto(3.195,0.56)(3.24,0.67)(3.28,0.74)
\curveto(3.32,0.81)(3.4,0.94)(3.44,1.0)
\curveto(3.48,1.06)(3.535,1.145)(3.55,1.17)
\curveto(3.565,1.195)(3.61,1.25)(3.64,1.28)
\curveto(3.67,1.31)(3.725,1.355)(3.75,1.37)
\curveto(3.775,1.385)(3.85,1.4)(3.9,1.4)
\curveto(3.95,1.4)(4.055,1.4)(4.11,1.4)
\curveto(4.165,1.4)(4.26,1.405)(4.3,1.41)
\curveto(4.34,1.415)(4.415,1.42)(4.45,1.42)
\curveto(4.485,1.42)(4.565,1.42)(4.61,1.42)
\curveto(4.655,1.42)(4.785,1.4)(4.87,1.38)
\curveto(4.955,1.36)(5.13,1.32)(5.22,1.3)
\curveto(5.31,1.28)(5.46,1.225)(5.52,1.19)
\curveto(5.58,1.155)(5.695,1.085)(5.75,1.05)
\curveto(5.805,1.015)(5.92,0.945)(5.98,0.91)
\curveto(6.04,0.875)(6.145,0.81)(6.19,0.78)
\curveto(6.235,0.75)(6.335,0.665)(6.39,0.61)
\curveto(6.445,0.555)(6.56,0.415)(6.62,0.33)
\curveto(6.68,0.245)(6.76,0.105)(6.78,0.05)
\curveto(6.8,-0.0050)(6.83,-0.125)(6.84,-0.19)
\curveto(6.85,-0.255)(6.86,-0.37)(6.86,-0.42)
\curveto(6.86,-0.47)(6.86,-0.57)(6.86,-0.62)
\curveto(6.86,-0.67)(6.845,-0.745)(6.83,-0.77)
\curveto(6.815,-0.795)(6.785,-0.84)(6.77,-0.86)
\curveto(6.755,-0.88)(6.725,-0.915)(6.71,-0.93)
\curveto(6.695,-0.945)(6.655,-0.985)(6.63,-1.01)
\curveto(6.605,-1.035)(6.55,-1.08)(6.52,-1.1)
\curveto(6.49,-1.12)(6.435,-1.16)(6.41,-1.18)
\curveto(6.385,-1.2)(6.325,-1.23)(6.29,-1.24)
\curveto(6.255,-1.25)(6.18,-1.265)(6.14,-1.27)
\curveto(6.1,-1.275)(6.035,-1.28)(6.01,-1.28)
\curveto(5.985,-1.28)(5.935,-1.285)(5.91,-1.29)
\curveto(5.885,-1.295)(5.82,-1.31)(5.78,-1.32)
\curveto(5.74,-1.33)(5.655,-1.345)(5.61,-1.35)
\curveto(5.565,-1.355)(5.48,-1.36)(5.44,-1.36)
\curveto(5.4,-1.36)(5.31,-1.37)(5.26,-1.38)
\curveto(5.21,-1.39)(5.075,-1.415)(4.99,-1.43)
\curveto(4.905,-1.445)(4.725,-1.475)(4.63,-1.49)
\curveto(4.535,-1.505)(4.36,-1.545)(4.28,-1.57)
\curveto(4.2,-1.595)(4.07,-1.63)(4.02,-1.64)
\curveto(3.97,-1.65)(3.885,-1.675)(3.85,-1.69)
\curveto(3.815,-1.705)(3.74,-1.74)(3.7,-1.76)
\curveto(3.66,-1.78)(3.57,-1.825)(3.52,-1.85)
\curveto(3.47,-1.875)(3.385,-1.91)(3.35,-1.92)
\curveto(3.315,-1.93)(3.255,-1.95)(3.23,-1.96)
\curveto(3.205,-1.97)(3.155,-1.995)(3.13,-2.01)
}
\rput(6.96,1.085){$N^*$}
\psline[linewidth=0.04cm,arrowsize=0.05291667cm 2.0,arrowlength=1.4,arrowinset=0.4]{->}(5.0,1.74)(6.18,-0.18)
\rput(4.8,0.8){$g\circ f$}
\end{pspicture} 
}
\end{center}  
\edem[Prop. \ref{EquiDomination}]


\begin{remark}
Notice that given $(M,\C{M},\tuple{a},N)$, if $M'$ is limit model over $M$ such that $N\indep_{M}^{\C{M}} M'$, in particular we have that $a\indep_{M}^{\C{M}}M'$ because $a\in N$. Therefore, if 
$\tuple{a}\dominates_M^{\C{M}} N$ we may say that $\tuple{a}$ and $N$ are equidominant over $M$ relative to $\C{M}$, which we denote by $\tuple{a}\equidominated_M^{\C{M}} N$.
\end{remark}

\begin{corollary}\label{Domination}
Given $(M,\C{M},\tuple{a},N)$ such that $\tuple{a}\indep_{M_\alpha}^{\C{M}_\alpha} M$ for some limit ordinal $\alpha$ such that $M_\alpha\in \C{M}$ (and therefore $\gatp(\tuple{a}/M)$ is a stationary type because $M$ is an universal model over $M_\alpha$), there exist $N^*$ and a resolution $\C{M}^*$ which witnesses that $M$ is a limit model over $M_0$ such that $\tuple{a}\equidominated_{M}^{\C{M}^*}N^*$.
\end{corollary}

\begin{question}\label{Question_Primes}
In general, we cannot assure the existence of prime models in metric
and discrete AECs. In superstable first order theories, we can prove
that if $p$ is a stationary syntactic type, there exist regular types
$p_1,\cdots , p_n$ such that $p\equidominated p_1\otimes \cdots
\otimes p_n$.   Setting $(\tuple{a}_1,\cdots,\tuple{a}_n)\models
p_1\otimes\cdots \otimes p_n$ and $\tuple{a}\models p$, it is known
that $M[\tuple{a}_1,\cdots \tuple{a}_n]=M[\tuple{a}]$ (i.e., the
a-prime model over $M\cup \{\tuple{a}\}$ and the a-prime model over
$M\cup\{\tuple{a}_1,\cdots,\tuple{a}_n\}$ agree). In Hilbert spaces
with a unitary operator (see~\cite{ArBe}), $a\indep_M N$ iff
$P_M(a)=P_N(a)$ (i.e., the respective orthogonal projections of $a$
over $M$ and $N$ agree). Considering this independence notion instead
of smooth-independence,  corollary~\ref{Domination} would say that
that given $\tuple{a}\in \C{H}$ (where $\C{H}$ is a monster Hilbert
space with a unitary operator) and $M$ a Hilbert space with a unitary
operator such that $M$ is saturated enough and $\tuple{a}\notin M$,
there exists a Hilbert space with a unitary operator $N^*\supset
acl(M\tuple{a})$ extending $M$ such that for every  Hilbert space with
a unitary operator $M'\ge M$, $P_M(\tuple{a})=P_{M'}(\tuple{a})$
implies that $P_M(b)=P_{M'}(b)$ for every $b\in N^*$; i.e.:
$\tuple{a}$ determines the projections on $M$ of all elements in
$N^*$.
A natural question that arises at this point is naturally connected to
the question on existence of prime models over sets in MAECs: under
which
assumptions can we guarantee that existence.


\end{question}

\subsection{Orthogonality}
Orthogonality arose from the question on the existence of bases
-maximal Morley sequences- of arbitrary size in a model, for (first
order syntactical types) $p$ and $q$ (see \cite{Bu}).

\indent In this section, we adapt the study of orthogonality which
S. Shelah did in the setting of good frames in (discrete) Abstract
Elementary Classes (see \cite{Sh705,Sh600}). Shelah provided a
suitable study of superstability in (discrete) AECs via good
frames, without assuming the   existence of a monster model as
in~\ref{Monster_Model} and with  an abstract notion of
independence. 
Most of the definitions in this section are inspired on Shelah's work
(\cite{Sh705}), with some exceptions (e.g., the definition of
domination of types, which we define in this thesis in order to prove
that domination corresponds to a kind of nonorthogonality). However,
we point out some differences between our results and the  analysis
done in \cite{Sh705}: although we are assuming the existence of a
homogeneous monster model (thereby losing generality), we are using a
fixed notion of
independence (smooth independence). In this section, we obtain an
adaptation of the notions given by Shelah to our setting and
prove some basic facts which were not proved in \cite{Sh705}.

However, we have to point out that there might be problems
proving the existence of weakly orthogonal types, with our definition
of weak orthogonality is being defined. Still, we
develop this section and show some important properties, consequence
of uniqueness of limit models.
%

\begin{notation}
$(M,\C{M},N,b,\alpha)$ means $\C{M}:=\{M_i:i<\delta\}$ is a resolution
of $M$ which witnesses that $M$ is a limit model, $\alpha<\delta$ is a limit ordinal, $N\ordenceac M$ is universal over $M$, $b\in N\setminus M$ and $b\indep^{\C{M}_\alpha}_{M_\alpha} M$, where $\C{M}_\alpha$ is a resolution of $M_\alpha$ such that $\C{M}_\alpha\subset \C{M}$.
\end{notation}

\begin{center}
\scalebox{0.6} 
{
\begin{pspicture}(0,-2.33)(6.8890624,2.33)
\psframe[linewidth=0.04,dimen=outer](6.3140626,2.33)(0.5940625,-2.33)
\psline[linewidth=0.04cm](3.3540626,2.27)(3.3740625,-2.27)
\psline[linewidth=0.04cm](0.6140625,-1.79)(3.3740625,-1.81)
\psline[linewidth=0.04cm](0.6340625,-1.17)(3.3540626,-1.17)
\psline[linewidth=0.04cm](0.6540625,-0.47)(3.3740625,-0.47)
\psline[linewidth=0.04cm](0.6340625,0.27)(3.3340626,0.31)
\psline[linewidth=0.04cm](0.6340625,0.93)(3.3340626,0.93)
\psline[linewidth=0.04cm](0.5940625,1.45)(3.3740625,1.49)
\psdots[dotsize=0.12](4.8140626,0.25)
\rput(4.6982813,-0.12){$b$}
\rput(6.723906,2.16){$N$}
\rput(0.14734375,1.44){$M$}
\rput(1.62,-0.04){$M_\alpha$}
\rput(1.62,-1.98){$M_0$}
\end{pspicture}
}
\end{center}

\subsubsection{Orthogonality and Independence of sequences}

\begin{definition}\label{Independence_Sequences}
Let $\mathbb{J}$ be a sequence of elements in $\monster$, $M\subset N$
where $M\ordencea N$ and $\C{M}$ be a resolution of $M$. We say that
{\it $\mathbb{J}$ is independent}
in $(M,N)$ iff there exist $\langle N_j,a_i : j\le \alpha,
i<\alpha \rangle$ such that
\begin{enumerate}
    \item $\langle N_i: i\le \alpha \rangle$ is a
      $\ordencea$-increasing and continuous chain.
    \item $\mathbb{J}=\{a_i:i<\alpha\}$.
    \item $M\ordencea N_i$
    \item $a_i\in N_{i+1}\setminus N_i$.
    \item $a_i\indep_{M}^{\C{M}} N_i$.
\end{enumerate}
\end{definition}

\begin{definition}
Let $M\in \C{K}$ be a limit model witnessed by $\C{M}$ and
$\C{M}_\alpha$ be a resolution of a model $M_\alpha\in \C{M}$ such
that $\C{M}_\alpha\subset \C{M}$.  Let $p,q\in \gaS(M)$ be
non-algebraic types such that $p,q\indep^{\C{M}_\alpha}_{M_\alpha}
M$. We say that {\it $p$ is {\it weakly orthogonal} to $q$ relative to
  $\alpha$} (denoted by $p\perp_\alpha^{wk} q$) iff
  for every $b\models q$ there exists $(M,\C{M},N,b,\alpha)$ where
$b\models q$ and if $p'\in \gaS(N)$ is any extension of $p$ then
$p'\indep_{M}^{\C{M}} {N}$ -notice that by definition of
$(M,\C{M},N,b,\alpha)$, $b\in N$-. We drop the subindex $\alpha$ if it
is clear.
\end{definition}

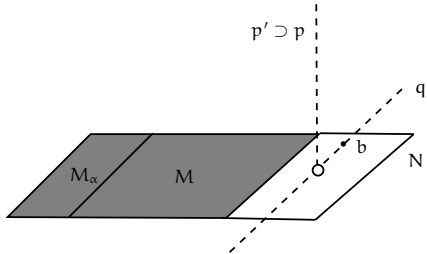
\scalebox{0.6} 
{
\begin{pspicture}(0,-2.77)(9.46,2.77)
\definecolor{color58b}{rgb}{0.6,1.0,0.2}
\definecolor{color58}{rgb}{1.0,0.0,0.2}
\pspolygon[linewidth=0.04,fillstyle=solid,fillcolor=gray](1.84,-0.13)(0.0,-1.97)(4.84,-1.99)(6.9,-0.13)
\psline[linewidth=0.04cm,fillcolor=color58b](6.92,-0.11)(9.0,-0.13)
\psline[linewidth=0.04cm,fillcolor=color58b](9.0,-0.13)(6.82,-2.03)
\psline[linewidth=0.04cm,fillcolor=color58b](6.82,-2.03)(4.76,-1.99)
\psline[linewidth=0.04cm,linecolor=black,fillcolor=color58b,linestyle=dashed,dash=0.16cm 0.16cm](4.92,-2.75)(8.74,0.87)
\psline[linewidth=0.04cm,linecolor=black,fillcolor=color58b,linestyle=dashed,dash=0.16cm 0.16cm](6.86,-0.85)(6.84,2.75)
\psline[linewidth=0.04cm,fillcolor=color58b](1.36,-1.97)(3.22,-0.11)
\usefont{T1}{ptm}{m}{n}
\rput(6,2.155){$p'\supset p$}
\psdots[dotsize=0.12](7.44,-0.35)
\usefont{T1}{ptm}{m}{n}
\rput(7.85,-0.405){$b$}
\usefont{T1}{ptm}{m}{n}
\rput(9.05,-0.705){$N$}
\usefont{T1}{ptm}{m}{n}
\rput(3.92,-1.065){$M$}
\usefont{T1}{ptm}{m}{n}
\rput(1.7,-1.045){$M_\alpha$}
\pscircle[linewidth=0.04,dimen=outer,fillstyle=solid](6.88,-0.93){0.14}
\usefont{T1}{ptm}{m}{n}
\rput(9.15,0.835){$q$}
\end{pspicture} 
}

In first order, for stationary types $p,q\in S(A)$, we say that $p$ is (almost) orthogonal to $q$ (over $A$) if and only if there exist realizations $a\models p$ and $b\models q$ such that $a\indep_A b$. Since in our setting we cannot consider independence either from or over sets which are not models, we have to adapt this notion to our context, as in \cite{Sh705}. Notice that if we could define smooth independence on sets, $p'\indep_M^{\C{M}} N$ would imply $p\indep_M^{\C{M}} b$ because $b\in N$. Hence, weak orthogonality corresponds to a stronger notion of orthogonality.   In spite of that, the notion of independence defined at the beginning of this subsection (independence of sequences, definition~\ref{Independence_Sequences}) allows us to catch such independence between a realization $a$ of $p$ and a realization $b$ of $q$.

\begin{example}\label{Ejemplo_Ortogonalidad}
Consider the class of Hilbert spaces. As in Hilbert spaces together with a unitary operator (see~\cite{ArBe}), independence is characterized by agreeing with the respective projections (i.e., $a\indep_N M$ if and only if $P_N(a)=P_M(a)$). In this case, replacing this notion of independence instead of smooth-independence, $p\perp^{wk} q$ (both of them in $\gaS(M)$) means that for \mbox{every} Hilbert space $N\ge M$ 
which contains a realization of $q$ and given any realization $a$ of $p$, $P_M(a)=P_N(a)$. If $M=\langle 0 \rangle$ and $N=\langle b \rangle$, notice that weak orthogonality implies that $0=P_M(a)=P_N(a)$, therefore $a$ and $b$ are orthogonal in the sense of the inner product in Hilbert spaces.  
\end{example}

\begin{proposition}\label{Orth_Indep}
Let $p,q\in\gaS(M)$ be non-algebraic types such that $p\perp_\alpha^{wk} q$ (witnessed by $(M,\C{M},N,b,\alpha)$)) and $\C{M}_\alpha$  and $\C{M}$ be resolutions of $M_\alpha$ and $M$ respectively such that $\C{M}_\alpha\subset \C{M}$, where $\C{M}$ witnesses that $M$ is a limit model. Therefore, for every $N'\ordenceac M$ and every realization $a$ of $p$ and a realization $b$ of $q$ belonging to $N'$ we have that $\langle b,a\rangle$ is independent in $(M,N')$.
\end{proposition}
\bdem
Define $N_0:=M$. Given $M$, $N'$, and $a,b\in N'$ as above, define $N_1:=N$. Notice that $b\notin N_0$ (since $q$ is non-algebraic). Trivially we have that $b\indep_{N_0}^{\C{N}_0} N_0$. Since $\gatp(a/N)\supset p$ and $(M,\C{M},N,b,\alpha)$ witnesses $p\perp_\alpha^{wk} q$  we may say $a\indep_{N_0}^{\C{N}_0} N$; i.e.: $a\indep_{M}^{\C{M}_0} N_1$. Notice that $a\notin N_1$ since $a\notin N_0$ (since $p$ is non-algebraic and by antirreflexivity, proposition \ref{s_antireflexivity}). Let $N_2\ordenceac N_1\cup N'$.
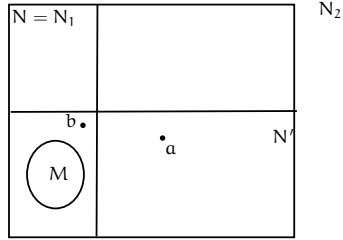
\begin{center}
\scalebox{0.6} 
{
\begin{pspicture}(0,-2.6091993)(7.6618943,2.6491992)
\psframe[linewidth=0.04,dimen=outer](6.36,2.5908008)(0.0,-2.6091993)
\psline[linewidth=0.04cm](1.98,2.5708008)(1.96,-2.5691993)
\psellipse[linewidth=0.04,dimen=outer](1.05,-1.1991992)(0.65,0.77)
\usefont{T1}{ptm}{m}{n}
\rput(1.1214551,-1.1641992){$M$}
\usefont{T1}{ptm}{m}{n}
\rput(0.7914551,2.2758007){$N=N_1$}
\psline[linewidth=0.04cm](0.06,0.19080079)(6.4,0.21080078)
\usefont{T1}{ptm}{m}{n}
\rput(6.101455,-0.3441992){$N'$}
\psdots[dotsize=0.12](1.66,-0.10919922)
\usefont{T1}{ptm}{m}{n}
\rput(1.391455,-0.004199219){$b$}
\psdots[dotsize=0.12](3.42,-0.38919923)
\usefont{T1}{ptm}{m}{n}
\rput(3.601455,-0.6441992){$a$}
\usefont{T1}{ptm}{m}{n}
\rput(7.131455,2.4558008){$N_2$}
\end{pspicture} 
}
\end{center}
Defining $a_0:=b$ and $a_1:=a$, notice that $\{N_0,N_1,N_2;a_0,a_1\}$ witness that $\langle b,a\rangle$ is independent in $(M,N')$.
\edem[Prop. \ref{Orth_Indep}]

The following proposition says that given $p,q\in \gaS(M)$) and $N\ordenceac M$ has a realization of $q$, then  $p$ is weakly orthogonal to $q$ if and only if $p$ has just one extension in $\gaS(N)$.

\begin{proposition}\label{UniqPerp}
Let $p,q\in \gaS(M)$ be non-algebraic types, $\C{M}$ be a resolution of $M$ which witnesses that $M$ is a limit model such that $p,q\indep^{\C{M}_\alpha}_{M_\alpha} M$ where  $M_\alpha\in \C{M}$ and $\C{M}_\alpha\subset \C{M}$ is a resolution of $M_\alpha$. Then $p\perp_\alpha^{wk} q$ $\Leftrightarrow$ for some $(M,\C{M},N,b,\alpha)$ such that $q=\gatp(b/M)$, $p$ has just one extension in $\gaS(N)$.
\end{proposition}
\bdem
\begin{enumerate}
\item[$\Rightarrow$] Suppose $p\perp_\alpha^{wk} q $ witnessed by $(M,\C{M},b,N,\alpha)$, then every extension $p'\in \gaS(N)$ of $p$ satisfies $p'\indep_{M}^{\C{M}} N$; since $M_\alpha$ and $M$ are limit model over $M_0$, by transitivity (proposition~\ref{transitivity})  we have $p'\indep_{M_\alpha}^{\C{M}_\alpha}N$. By stationarity over limit models of s-independence (proposition \ref{UnExtUniv}), we have that there is just one extension of $p$ in $\gaS(N)$.
\item[$\Leftarrow$] Let $(M,\C{M},b,N,\alpha)$ be such that $b\models q$ and such that $p$ has just one extension in $\gaS(N)$.  By extension property of s-independence (proposition \ref{ExtUniv}, since $M$ is a limit model over $M_\alpha$), there exists an extension $p'\supset p$ in $\gaS(N)$ such that $p'\indep_{M_\alpha}^{\C{M}_\alpha} N$. By monotonicity (proposition \ref{monotonicity} , since $\C{M}_\alpha\subset \C{M}$), $p'\indep_{M}^{\C{M}}N$. Since $p$ has just one extension in $\gaS(N)$, then the unique extension of $p$ in $\gaS(N)$  is $p'$. Therefore,  $(M,\C{M},b,N,\alpha)$ witnesses $p\perp_\alpha^{wk} q$.
\end{enumerate}
\edem[Prop. \ref{UniqPerp}]

\begin{proposition}\label{Not_perp}
Let $(M,\C{M},N,b,\alpha)$ be such that $b\models q$ and $p\in \gaS(M)$ be a non-algebraic Galois type such that $p\indep_{M_\alpha}^{\C{M}_\alpha} M$. If $p$ is realized in $N$, then  $(M,\C{M},N,b,\alpha)$ cannot witness $p\perp_\alpha^{wk} q$.
\end{proposition}
\bdem
Since $M$ is universal over $M_\alpha$ and $p\indep_{M_\alpha}^{\C{M}_\alpha} M$, by extension property (proposition~\ref{ExtUniv}) there exists $p'\supset p$ in $\gaS(N)$ such that $p'\indep_{M_\alpha}^{\C{M}_\alpha} N$. Notice that $p'$ is non-algebraic: otherwise, by antirreflexivity  (proposition~\ref{s_antireflexivity} (6), since $M_\alpha$ is a limit model witnessed by $\C{M}_\alpha$) $p'$ would be realized in $M_\alpha$ and so realized in $M$ (contradiction). But by hypothesis, there exists $c\in N$ such that $c\models p$. Notice that $p'':=\gatp(c/N)\supset p$ and $p'\neq p''$. If $(M,\C{M},N,b,\alpha)$ witnessed $p\perp_\alpha^{wk} q$, this would contradict~proposition~\ref{UniqPerp}.
\edem[Prop. \ref{Not_perp}]

Another way to understand the previous proposition is the following\linebreak
\mbox{corollary}:

\begin{corollary}  
Let $(M,\C{M},N,b,\alpha)$ be such that $b\models q$ and $p\in \gaS(M)$ be a non-algebraic Galois type such that $p\indep_{M_\alpha}^{\C{M}_\alpha} M$. If  $p\perp_\alpha^{wk} q$ is witnessed by $(M,\C{M},N,b,\alpha)$, then $p$ is not realized in $N$.
\end{corollary}

As we stated in section~\ref{section:domination}, orthogonality corresponds (in first order) in some way to nonorthogonality. In order to prove a similar result in our context, we adapt the notion of domination of types.  

\begin{definition}
Let $p,q\in \gaS(M)$ be non-algebraic Galois types such that $p,q\indep_{M_\alpha}^{\C{M}_\alpha} M$. We say that  $q$ {\it is dominated by} $p$ (denoted by $q\isdominated p$) if there exist $a\models p$, $b\models q$ and $N\ordenceac M$ such that $a,b\in N$ and $N\isdominated_M a$.
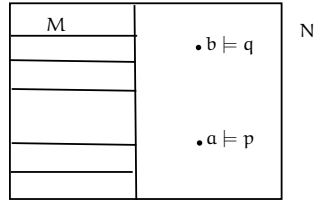
\begin{center}
\scalebox{0.6} 
{
\begin{pspicture}(0,-2.21)(6.96,2.19)
\psframe[linewidth=0.04,dimen=outer](6.06,2.19)(0.0,-2.19)
\psline[linewidth=0.04cm](2.78,2.13)(2.82,-2.19)
\psdots[dotsize=0.12](4.2,1.19)
\psdots[dotsize=0.12](4.22,-0.93)
\rput(4.9,1.235){$b\models q$}
\rput(4.9,-0.845){$a\models p$}
\rput(6.61,1.575){$N$}
\rput(1.04,1.7){$M$}
\psline[linewidth=0.04cm](0.04,1.45)(2.84,1.45)
\psline[linewidth=0.04cm](0.02,0.91)(2.8,0.87)
\psline[linewidth=0.04cm](0.06,0.27)(2.82,0.23)
\psline[linewidth=0.04cm](0.04,-1.57)(2.74,-1.57)
\psline[linewidth=0.04cm](0.06,-0.93)(2.82,-0.97)
\end{pspicture} 
}
\end{center}

\end{definition}

The following propositions says that this notion of domination is a kind of ``opposite'' of weak orthogonality.

\begin{proposition}\label{Orthogonal_Domination}
Let  $\C{M}:=\{M_i:i<\theta\}$ be a resolution of a model $M$ and $\C{M}_\alpha\subset \C{M}$ be a resolution of $M_\alpha$. If $q\isdominated p$ witnessed by $a\models p$, $b\models q$ and $a,b\in N$, then $(M,\C{M},N,b,\alpha)$ does not witness $p\perp_\alpha^{wk}q $ and  $(M,\C{M},N,a,\alpha)$ does not witness $q\perp^{wk}_{\alpha}p$.
\end{proposition}
\bdem
Let $a\models p$, $b\models q$ and $N\ordenceac M$ be witnesses of $q\isdominated p$. Notice that $p':=\gatp(a/N)$ is an extension of $p$ such that $a\not\hspace{-2.4mm}\indep^{\C{M}}_{M} N$ (by antirreflexivity -proposition~\ref{s_antireflexivity}-, since $a\in N\setminus M$). Notice that $(M,\C{M},N,b,\alpha)$ does not witness  $p\perp^{wk}_{\alpha} q$. By an analogous argument, we can prove that $(M,\C{M},N,a,\alpha)$ does not witness $q\perp^{wk}_{\alpha} p$, using $(M,\C{M},N,a,\alpha)$.
\edem[Prop. \ref{Orthogonal_Domination}]

\begin{remark}
Notice that proposition~\ref{Orthogonal_Domination} says that if given $a\models p$ and $b\models q$, if $(M,\C{M},N,b,\alpha)$ witnesses $p\perp^{wk} q$ and $(M,\C{M},N,a,\alpha)$ witnesses $q\perp^{wk} p$, then $a,b,N$ cannot witness $q\isdominated p$.
\end{remark}

%
%

\subsection{Parallelism}
Roughly speaking, two (first order syntactical) stationary types $p$ and $q$ are parallel if and only if they have a common independent extension. In this section, we study parallelism of strong limit Galois types in the setting of superstable MAECs.
\\ \\
\indent We defined parallelism of strong Galois types in \cite{Za,ViZa} inspired in the definition given in  \cite{GrVaVi}, which we used it as an auxiliary tool for studying full-relative s-towers. Full-relative s-towers were very important to get a proof of uniqueness of limit models because they codified a kind of saturation. In this subsection, we study some properties of a stronger version of parallelism, but in the setting of superstable MAECs.
\\ \\
\indent However, for the sake of completeness, we provide the definition of parallelism once more. But we have to point out that in this subsection, we require that if $(p,N)\in \stype(M)$, then $M$ is a limit model over $N$, instead of just being a universal model over $N$. Because of that, we define a stronger notion of strong type, which we call {\it strong limit type}. In this thesis, we use the notion of {\it strong limit type}
 instead of {\it strong types} because we want to use {\it uniqueness of limit models}  to prove some properties of parallelism, e.g. proposition \ref{Prop_Orth_Parall} (2) and (3).

\begin{definition}[strong limit type]\label{StrongLimitType}
Let $M$ be a $(\mu,\sigma)$-limit model\\
$$\mathfrak{SL}(M):=\left\{
(p,N) :
\begin{tabular}{l}
$N\ordencea M$\\
$N$ is a $\theta$-Limit Model\\
$M$ is a Limit Model over $N$\\
$p\in \gaS(M)$ is non-algebraic\\
and $p\indep^{\C{N}}_{N} M $\\
for some resolution $\C{N}$ of $N$.
\end{tabular}\right\}$$
\end{definition}

\begin{definition}[Parallelism]\label{Parallelism1}
Two strong limit  types $(p_l,N_l)\in \mathfrak{SL}(M_l)$ ($l\in
\{1,2\}$) are said to be {\it parallel}  (which we denote by
$(p_1,N_1)\parallel (p_2,N_2)$) iff for every $M'\ordenceac
{M}_1,{M}_2$ with density character $\mu$, there exists $q\in
\gaS(M')$ which extends both $p_1$ and $p_2$ and
$q\indep^{\C{N}_l}_{N_l} M'$ ($l\in\{1,2\}$)
(where $\C{N}_l$ is the resolution of $N_l$ which satisfies
$p_l\indep^{\C{N}_l}_{N_l} M_l$). If there is no any confusion, we denote it by $p_1\parallel p_2$.
\end{definition}  

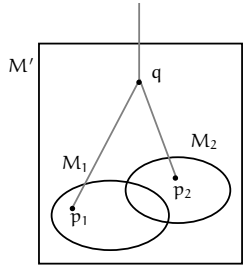
\begin{center}
\scalebox{0.6} 
{
\begin{pspicture}(0,-2.9)(4.9140625,2.92)
\psframe[linewidth=0.04,dimen=outer](4.9140625,2.02)(0.3740625,-2.9)
\psellipse[linewidth=0.04,dimen=outer](1.9740624,-1.81)(1.32,0.79)
\psellipse[linewidth=0.04,dimen=outer](3.4740624,-1.25)(1.18,0.75)
\rput(1.2,-0.67){$M_1$}
\rput(4.0554686,-0.17){$M_2$}
\rput(0,1.53){$M'$}
\psline[linewidth=0.04cm,linecolor=gray](1.1340625,-1.64)(2.6340625,1.22)
\psline[linewidth=0.04cm,linecolor=gray](3.4540625,-1.02)(2.6340625,1.18)
\psline[linewidth=0.04cm,linecolor=gray](2.6140625,1.18)(2.6140625,2.9)
\psdots[dotsize=0.12,linecolor=black](3.4140625,-0.98)
\psdots[dotsize=0.12,linecolor=black](1.1340625,-1.66)
\psdots[dotsize=0.12,linecolor=black](2.6140625,1.14)
\rput(1.2985938,-1.9){$p_1$}
\rput(3.578125,-1.3){$p_2$}
\rput(3,1.29){$q$}
\end{pspicture}
}
\end{center}

\begin{remark}
Consider the class of Hilbert spaces.  Let us suppose that we could set $N_1=N_2=\langle 0 \rangle \subset \mathbb{R}^3$ -the space generated by the origin- (despite this is not a universal model) and let $M:=M_1=M_2=\{(x,0,0): x\in \mathbb{R}\}$. Remember that we stated in~\ref{Ejemplo_Ortogonalidad} that, as in a  
Hilbert space with a unitary operator (see \cite{ArBe}), independence in Hilbert spaces means the respective projections agree.  Let $p_i\in \gaS(M)$, $a\models p_1$ and $b\models p_2$ be such that $a$ and $b$ are independent from $M$ over $\langle 0 \rangle$; i.e.: $0=P_M(a)=P_M(b)$, therefore $a$ and $b$ are orthogonal to $M$.

If $p_1$ and $p_2$ are parallel in the sense defined above, consider $M':=\langle M\cup\{a,b\} \rangle\ge M$, so there exists a type over $M'$ $q\supset p_1,p_2$ (so $\gatp(a/M)=\gatp(b)/M)=q\rest M$) such that $q$ is independent from $M'$ over $\langle 0\rangle$ (i.e., any realization $c\models q$ satisfies $0=P_{M'}(c)$), therefore $c$ is orthogonal to $a$ and $b$. If $c\in \mathbb{R}^3$, notice that it means that if $\varphi$ is the angle between $a$ and b and $\theta$ is the angle between $b$ and $c$ (and so $\theta+\varphi$ is the angle between $a$ and $c$),  since $a$ and $c$ are orthogonal then $|\cos(\varphi+\theta)|=0$, and since $b$ and $c$ are orthogonal then $|\cos(\theta)|=\cos(\theta)=0$, then $\theta=\frac{(2k+1)\pi}{2}$ for some $k\in \mathbb{Z}$. Since $0=|\cos(\varphi+\theta)|=|\cos(\theta)\cos(\varphi)-\sin(\theta)\sin(\varphi)|=|\sin\left(  \frac{(2k+1)\pi}{2} \right)\sin(\varphi)|=|\sin(\varphi)|$, therefore $\varphi=m\pi$ for some $m\in \mathbb{Z}$; i.e.: $a$ and b would be parallel as vectors in $\mathbb{R}^3$.

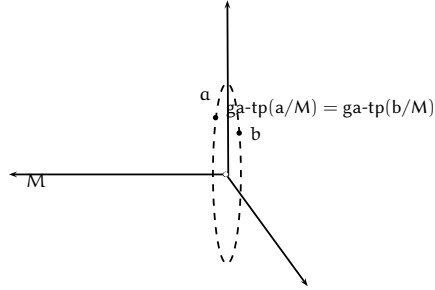
\begin{center}
\scalebox{0.6} 
{
\begin{pspicture}(0,-4.44)(11.3,4.44)
\psline[linewidth=0.04cm,arrowsize=0.05291667cm 2.0,arrowlength=1.4,arrowinset=0.4]{->}(4.86,0.52)(4.84,4.42)
\psline[linewidth=0.04cm,arrowsize=0.05291667cm 2.0,arrowlength=1.4,arrowinset=0.4]{->}(4.84,0.56)(0.0,0.56)
\psellipse[linewidth=0.04,linestyle=dashed,dash=0.16cm 0.16cm,dimen=outer](4.83,0.58)(0.33,2.02)
\psline[linewidth=0.04cm,arrowsize=0.05291667cm 2.0,arrowlength=1.4,arrowinset=0.4]{->}(4.84,0.52)(6.62,-1.92)
\psdots[dotsize=0.12,fillstyle=solid,dotstyle=o](4.8,0.56)
\psdots[dotsize=0.12,fillstyle=solid,dotstyle=o](11.22,-4.36)
\rput(7.13,2.025){$\gatp(a/M)=\gatp(b/M)$}
\rput(0.6,0.425){$M$}
\psdots[dotsize=0.12](4.58,1.82)
\psdots[dotsize=0.12](5.1,1.48)
\rput(4.34,2.305){$a$}
\rput(5.47,1.465){$b$}
\end{pspicture} 
}
\end{center}

\begin{center}
\scalebox{0.6} 
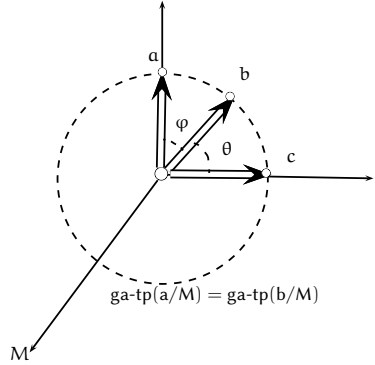
{
\begin{pspicture}(0,-3.97)(7.86,3.95)
\pscircle[linewidth=0.04,linestyle=dashed,dash=0.16cm 0.16cm,dimen=outer](3.17,-0.02){2.35}
\psdots[dotsize=0.3,fillstyle=solid,dotstyle=o](3.14,0.09)
\psline[linewidth=0.04cm,arrowsize=0.05291667cm 2.0,arrowlength=1.4,arrowinset=0.4]{->}(3.06,-0.05)(0.22,-3.87)
\psline[linewidth=0.04cm,arrowsize=0.05291667cm 2.0,arrowlength=1.4,arrowinset=0.4]{->}(3.28,0.07)(7.84,-0.01)
\psline[linewidth=0.04cm,arrowsize=0.05291667cm 2.0,arrowlength=1.4,arrowinset=0.4]{->}(3.16,0.23)(3.16,3.93)
\psline[linewidth=0.04cm,arrowsize=0.05291667cm 2.0,arrowlength=1.4,arrowinset=0.4,doubleline=true,doublesep=0.12]{->}(3.12,0.23)(3.16,2.35)
\psline[linewidth=0.04cm,arrowsize=0.05291667cm 2.0,arrowlength=1.4,arrowinset=0.4,doubleline=true,doublesep=0.12]{->}(3.28,0.17)(4.76,1.81)
\psdots[dotsize=0.2,fillstyle=solid,dotstyle=o](3.16,2.35)
\psdots[dotsize=0.2,fillstyle=solid,dotstyle=o](4.66,1.81)
\rput(4.33,-2.605){$\gatp(a/M)=\gatp(b/M)$}
\rput(0,-3.9){$M$}
\rput(2.96,2.695){$a$}
\rput(4.99,2.315){$b$}
\psline[linewidth=0.04cm,arrowsize=0.05291667cm 2.0,arrowlength=1.4,arrowinset=0.4,doubleline=true,doublesep=0.12]{->}(3.34,0.09)(5.58,0.07)
\pscustom[linewidth=0.04,linestyle=dashed,dash=0.16cm 0.16cm]
{
\newpath
\moveto(3.18,0.85)
\lineto(3.23,0.85)
\curveto(3.255,0.85)(3.3,0.84)(3.32,0.83)
\curveto(3.34,0.82)(3.38,0.795)(3.4,0.78)
\curveto(3.42,0.765)(3.46,0.745)(3.48,0.74)
\curveto(3.5,0.735)(3.535,0.715)(3.55,0.7)
\curveto(3.565,0.685)(3.59,0.655)(3.62,0.61)
}
\psdots[dotsize=0.2,fillstyle=solid,dotstyle=o](5.46,0.11)
\rput(6.02,0.455){$c$}
\rput(3.6,1.215){$\varphi$}
\pscustom[linewidth=0.04,linestyle=dashed,dash=0.16cm 0.16cm]
{
\newpath
\moveto(3.86,0.77)
\lineto(3.9,0.76)
\curveto(3.92,0.755)(3.955,0.735)(3.97,0.72)
\curveto(3.985,0.705)(4.015,0.675)(4.03,0.66)
\curveto(4.045,0.645)(4.075,0.61)(4.09,0.59)
\curveto(4.105,0.57)(4.135,0.53)(4.15,0.51)
\curveto(4.165,0.49)(4.18,0.445)(4.18,0.42)
\curveto(4.18,0.395)(4.185,0.34)(4.19,0.31)
\curveto(4.195,0.28)(4.2,0.22)(4.2,0.19)
\curveto(4.2,0.16)(4.2,0.125)(4.2,0.11)
}
\rput(4.6,0.655){$\theta$}
\end{pspicture} 
}
\end{center}
\end{remark}

\begin{claim}\label{Equiva_parallel}
$\parallel$ is an equivalence relation.
\end{claim}

\bdem
Reflexivity and symmetry are trivial. We focus on transitivity. Let $(p_l,N_l)\in  \mathfrak{SL}(M_l)$ ($l\in \{1,2,3\}$) 
be such that  $p_1\parallel p_2$ and $p_2\parallel p_3$. Let $M$ be a $\C{K}$-extension of both $M_1$ and $M_3$. By Downward L\"owenheim-Skolem axiom (Definition  \ref{MAEC} (6)) and Coherence axiom of MAEC (Definition \ref{MAEC} (5)), there exists a $\C{K}$-extension $M'$ of both $M$ and $M_2$. Denote by $p^{M'}_l$ ($l\in \{1,2,3\}$) the unique s-independent extension of $p_l$ in $\gaS(M')$ (such that extension 
\mbox{exists} by propositions~\ref{ExtUniv} and \ref{UnExtUniv}, since $M$ is universal over $N_l$) and by $p_k^{M}$ ($k\in \{1,3\}$) the unique s-independent extension of $p_k$ in $\gaS(M)$. Notice that for $k\in \{1,3\}$ we have $p_k^{M}=p_k^{M'}\rest M$. Since $M_1,M_2\ordencea M'$ and $p_1\parallel p_2$, then $p_1^{M'}=p_2^{M'}$. In a similar way we have $p_2^{M'}=p_3^{M'}$, then $p_1^{M'}=p_3^{M'}$. Since $M\ordencea M'$, then $p_1^{M}=p_1^{M'}\rest M = p_3^{M'}\rest M =p_3^M$. Therefore, $p_1\parallel p_3$. \ \ \ \
\edem[Claim \ref{Equiva_parallel}]

The following proposition says that strong limit types are stationary (up to parallelism).
\begin{proposition}[``Stationarity'' of parallelism]\label{Stat_parallel}
Let $(p,N)\in \mathfrak{SL}(M)$ and $M'\ordenceac M$ be a limit model over $M$. There exists a unique $(q,N)\in \mathfrak{SL}(M')$ such that $p\parallel q$.%
\end{proposition}
\bdem
Since $M$ is universal over $N$, by stationarity (proposition~\ref{UnExtUniv}) there exists a unique $q\in \gaS(M')$ such that $q\indep_N^{\C{N}} M'$. Notice that $(q,N)\in \mathfrak{SL}(M')$.

We have that $p\parallel q$: Let $M''\ordenceac M'$ (and so $M''\ordenceac M$). If $p',q'\in \gaS(M')$ are the s-independent extensions of $p$ and $q$ respectively, we have $p'=q'$ (if not, $p'\neq q'$ are s-independent extensions of $p$, contradicts stationarity). Therefore $p\parallel q$.

If $(q^*,N)\in \mathfrak{SL}(M')$ satisfies $p\parallel q^*$ and $q`\in \gaS(M')$ is the unique extension of $p$ and $q^*$ (so $q^*=q'$) such that $q'\indep_N^{\C{N}} M'$, then by stationarity (proposition~\ref{UnExtUniv}) $q^*=q'=q$. So, uniqueness is proved.
\edem[Prop. \ref{Stat_parallel}]

Next, we prove that weak orthogonality is preserved under parallelism. Before giving its proof, we prove that weak orthogonality is invariant under isomorphisms and that weak orthogonality is preserved under $\C{K}$-submodels and $\C{K}$-superstructures if we have suitable independence conditions.

\begin{proposition}\label{Prop_Orth_Parall}
\begin{enumerate}
\item Given $p, q \in\gaS(M)$, $\C{M}:=\{M_i:i<\delta\}$ a resolution of $M$ which witnesses that $M$ is a limit model such that $p,q\indep_{M_\alpha}^{\C{M}_\alpha} M$ for some $\alpha<\delta$ (where $\C{M}_\alpha\subset \C{M}$ is a resolution of $M_\alpha$) and $f:M\approx N$ is an isomorphism, then $(M,\C{M},N',b,\alpha)$ witnesses $p\perp_\alpha^{wk} q$ if and only if  $(N,f[\C{M}],\overline{f}[N'],\overline{f}(b),\alpha)$ witnesses $ f(p)\perp_\alpha^{wk} f(q)$, where $\overline{f}\in Aut(\monster)$ extends $f$.
\item Given $\C{M}:=\{M_i:i<\delta\}$ a resolution which witnesses that $M$ is a limit model,  if $N\ordenceac M$ is limit over $M$, given $p,q\in \gaS(N)$ such that $p,q\indep_{M_\alpha}^{\C{M}_\alpha} N$, $p\perp_\alpha^{wk} q \Leftrightarrow p\rest M \perp_\alpha^{wk} q\rest M$.
\item Given $\C{M}:=\{M_i:i<\delta\}$ a resolution which witnesses that $M$ is a limit model, if $N\ordenceac M$ is a limit model over $M$ (and in particular over $M_{\alpha+1}$) and $(p_l,M_{\alpha})\in \mathfrak{SL}(M)$ and $(q_l,M_{\alpha})\in \mathfrak{SL}(N)$ ($l\in\{1,2\}$) satisfy $p_i\parallel q_i$ ($i\in \{1,2\}$), then $p_1\perp_\alpha^{wk} p_2$ iff $q_1\perp_\alpha^{wk} q_2$.
\end{enumerate}
\end{proposition}

\bdem
\begin{enumerate}
\item Since $(M,\C{M},N,b,\alpha)$ witnesses $p\perp^{wk} q$ then $p\indep_{M_\alpha}^{\C{M}_\alpha} N$. Notice that $f[\C{M}]:=\{f[M_i]:i<\delta\}$ witnesses that $N$ is a limit model. By invariance (fact~\ref{invariance}), $f(p)\indep_{f[M_\alpha]}^{f[\C{M}]_\alpha} f[N]$, therefore\linebreak
$(N,f[\C{M}],f[N'],f(b),\alpha)$ witnesses $f(p)\perp_\alpha^{wk} f(q)$. Converse follows from a similar argument as above.

\item Let $\C{N}^*$ be a resolution of $N$ witnessing $N$ is a limit model over $M$, so $\C{M}^\frown \C{N}^*$ witnesses that $N$ is a limit model over $M_0$ (and in particular, it witnesses that $N$ is a limit model over $M_{\alpha+1}$).  Since $M$ and $N$ are limit models over $M_{\alpha+1}$, by assumption~\ref{Uniqueness} (uniqueness of limit models) there exists $f:M\approx_{M_{\alpha+1}} N$.   Notice that $p\rest M_{\alpha+1}=f(p\rest M_{\alpha+1})\subset f(p\rest M)\in \gaS(N)$ and $q\rest M_{\alpha+1}=f(q\rest M_{\alpha+1})\subset f(q\rest M)\in \gaS(N)$. Since $p\rest M_{\alpha+1} \indep_{M_\alpha}^{\C{M}_\alpha} M_{\alpha+1}$ (by monotonicity, since $p\indep_{M_\alpha}^{\C{M}_\alpha} N$) and $f(p\rest M)\supset p\rest M_{\alpha+1}$ satisfies $f(p\rest M)\indep_{M_\alpha}^{\C{M}_\alpha} N$ (by invariance of s-independence, since $p\rest M\indep_{M_\alpha}^{\C{M_\alpha}} M$) and we also have $p\indep_{M_\alpha}^{\C{M}_\alpha} N$, then by stationarity of $p\rest M_{\alpha+1}$ (proposition \ref{UnExtUniv}, notice that $M_{\alpha+1}$ is universal over $M_\alpha$) we have that $f(p\rest M)=p$. In a similar way we can prove $f(q\rest M)=q$. By proposition \ref{Prop_Orth_Parall} (1) we have $(M,\C{M},N',b,\alpha)$ witnesses $p\rest M\perp_\alpha^{wk} q\rest M$ iff  $(N,f[\C{M}],\overline{f}[N'],\overline{f}(b),\alpha)$ witnesses $p\perp_\alpha^{wk} q$, where $\overline{f}\in Aut(\monster)$ extends $f$.

\item Notice that in this case, $p_i\parallel q_i$ implies $p_i= q_i\rest M$. So, this holds by proposition \ref{Prop_Orth_Parall} (2).
\end{enumerate}
\edem[Prop. \ref{Prop_Orth_Parall}]

\end{document}